\documentclass[a4paper,fleqn]{cas-sc}

\usepackage[numbers]{natbib}
\usepackage{amsmath}
\usepackage{amssymb}
\usepackage{multirow}
\usepackage{framed,color}
\definecolor{shadecolor}{rgb}{0.9,0.9,0.9}

\newcommand{\R}{\,{\mathbb R}}
\newcommand{\C}{\,{\mathbb C}}

\newcommand{\T}{\,{\mathbb T}}
\newcommand{\diag}{\mbox{diag}}

\newcommand{\calH}{{\mathcal H}}
\newcommand{\calP}{{\mathcal P}}

\renewcommand{\vec}{\mbox{vec}}

\newcommand{\vecf}{{\mathbf f}}
\newcommand{\vecF}{{\mathbf F}}
\newcommand{\vecg}{{\mathbf g}}
\newcommand{\vecv}{{\mathbf v}}
\newcommand{\vecx}{{\mathbf x}}
\newcommand{\xbar}{\overline{x}}
\newcommand{\thetabar}{\overline{\theta}}
\newcommand{\vecxbar}{\overline{\mathbf x}}
\newcommand{\veczero}{{\mathbf 0}}

\newenvironment{new}{\color{black}}{\color{black}}
 \newenvironment{moved}{\color{black}}{\color{black}}
 \newenvironment{reorganized}{\color{black}}{\color{black}}\date{}

\begin{document}
\shorttitle{Low Order Spectral Factorization}
\shortauthors{Kolev et~al.}

\title[mode=title]{Bauer's Spectral Factorization Method for Low Order
  Multiwavelet Filter Design}

% \author{Vasil Kolev\corref{cor1}}
\author{Vasil Kolev}
\cormark[1]
\ead{vasil.kolev@iict.bas.bg, kolev_acad@abv.bg}
\address{Institute of Information and Communication Technologies,
  Bulgarian Academy of Sciences, Bl. 2, Acad. G. Bonchev St., 1113
  Sofia, Bulgaria}
%\cortext{Corresponding author}

\author{Todor Cooklev}[orcid=0000-0002-4619-3579]
\address{Wireless Technology Center, Purdue University, Fort Wayne, IN 46805, USA}
\ead{cooklevt@pfw.edu}

\author{Fritz Keinert}[orcid=0000-0003-3116-967X]
\address{Department of Mathematics, Iowa State University, Ames, IA 50011, USA} 
\ead{keinert@iastate.edu}

\begin{abstract}
  Para-Hermitian polynomial matrices obtained by matrix spectral
  factorization lead to functions useful in control theory systems,
  basis functions in numerical methods or multiscaling functions used in
  signal processing. We introduce a fast algorithm for matrix spectral
  factorization based on Bauer’s method.  We convert Bauer’s method into
  a nonlinear matrix equation (NME). The NME is solved by two different
  numerical algorithms (Fixed Point Iteration and Newton’s Method) which
  produce approximate scalar or matrix factors, as well as a symbolic
  algorithm which produces exact factors in closed form for some
  low-order scalar or matrix polynomial matrices,
  respectively. Convergence rates of the two numerical algorithms are
  investigated for a number of singular and nonsingular scalar and
  matrix polynomials taken from different areas. In particular, one of
  the singular examples leads to new orthogonal multiscaling and
  multiwavelet filters. Since the NME can also be solved as a
  Generalized Discrete Time Algebraic Riccati Equation (GDARE),
  numerical results using built-in routines in Maple 17.0 and six Matlab
  versions are presented.
\end{abstract}


\maketitle

\section{Introduction and Problem Statement}\label{sec:intro}

{\em Spectral Factorization} (SF) is an important tool in signal
processing, control theory, information theory, and many other
fields. Despite the vast literature on the subject, there is still a lot
of interest in this topic.  SF factorizes a given Laurent polynomial
function of the form
\begin{equation*}
  p(z) = \sum_{k=-m}^n p_k z^k
\end{equation*}
into causal and anti-causal polynomials. As the name indicates, it is
often applied to spectral densities of stochastic processes.  In the
case of polynomials with matrix coefficients, SF is called {\em Matrix
  Spectral Factorization} (MSF).

In systems and control theory, para-Hermitian matrix polynomials play a
fundamental role. The para-conjugate of a matrix Laurent polynomial is
defined as
\begin{equation*}
  P(z)^* = \left( \sum_{k=-m}^n P_k z^k \right)^* =  \sum_{k=-m}^n P_k^* z^{-k},
\end{equation*}
where $A^*$ is the complex conjugate transpose of the matrix $A$.  A
matrix polynomial is {\em para-Hermitian} if $P(z)^* = P(z)$.  Thus, a
para-Hermitian matrix polynomial has the form
\begin{equation*}
  P(z) = \sum_{k=-m}^m P_k z^k, \quad P_{-k} = P_k^*.
\end{equation*}
Such matrix polynomials were first used by Norbert Wiener
\cite{W49,WM57}, where they represent spectral density functions of
stochastic processes.  They are also used in many other settings, for
example in MIMO communication theory, in linear quadratic optimal
control, optimal filtering engineering systems, signal processing, and
control analysis.

A matrix $A$ is {\em positive definite}, written $A > 0$, if
$\vecx^* A \vecx > 0$ for all $\vecx \ne \veczero$.  It is {\em positive
  semidefinite}, written $A \ge 0$, if $\vecx^* A \vecx \ge 0$ for all
$\vecx$.  The matrix form of the Riesz-Fej{\'e}r lemma
\cite{EJL09,HHS04} implies that if $P(z)$ is para-Hermitian and positive
semidefinite for $z$ on the unit circle $\T$ in the complex plane, then
it can be {\em left factorized} as $P(z) = H(z) H(z)^*$,
$z \in \C \backslash \{0\}$, where
\begin{equation*}
  H(z) = \sum_{k=0}^m H_k z^{-k}.
\end{equation*}
$H(z)$ is called a {\em left spectral factor} of $P(z)$.

Similarly, a {\em right spectral factor} satisfies $P(z) = H(z)^* H(z)$.
In this paper we consider left factorizations only.

A matrix polynomial $P(z)$ is called {\em singular} if its determinant
$|P(z)|$ vanishes somewhere on $\T$.  Singular polynomial matrices
appear for example in control and linear systems theory \cite{K80},
wireless communication \cite{BML03}, and especially in wavelet and
multiwavelet filter bank design \cite{KCK18}.

In multiwavelet theory, scalar or matrix spectral factorization plays a
key role.  Finding spectral factors of para-Hermitian polynomial
matrices (which are product filters in signal processing theory) leads
to desirable multiscaling functions satisfying simultaneously important
wavelet properties -- symmetry, smoothness, short support, and
approximation order.  There exist different spectral factorization
methods \cite{ESS18, YK78, WG69, BFGM03, RDF89, PVBH09}, but none leads
to exact values of spectral factors. Factorization is a key step in the
construction of wavelet or multiwavelet filter banks with high precision
coefficients. This approach was first used in Cooklev’s PhD Thesis
\cite{C95} for the construction of orthogonal multiwavelet filter banks.
In that application, the smoothness of the multiwavelet is improved if
the determinant of $P(z)$ has a higher-order zero at $z = -1$ (see
\cite{V86}).

Most of the numerical methods for MSF proposed in the past only work in
the nonsingular case, that is, the case where $|P(z)|$ does not vanish
on $\T$. Bauer's method, as presented in Youla and Kazanjian
\cite{YK78}, is based on the Cholesky decomposition of a certain
semi-infinite block Toeplitz matrix. It continues to work even in the
singular case. However, it converges very slowly in this case, and the
result has reduced accuracy.

Our goal in this paper is to present modifications of Bauer's method
that speed up convergence in the singular case, and even produce exact
solutions in some simple situations. Our work can be understood as an
extension of Youla and Kazanjian's theoretical work, providing new
results on convergence speed and symbolic solution. Both the numerical
and symbolic methods are based on solving a nonlinear matrix equation
(NME) introduced below.

We have the following main contributions:
\begin{itemize}
\item We construct algorithms for finding spectral factorizations of
  singular and nonsingular scalar and matrix polynomials; one algorithm
  produces exact values of spectral factors for low order scalar and
  matrix polynomials, and two other algorithms produce numerical
  approximations (based on fixed point iteration (FPI) or Newton’s
  method);
\item We consider for the first time in detail the convergence speed and
  accuracy in the numerical algorithms;
\item We obtain a new orthogonal multiwavelet filter from factoring the
  para-Hermitian polynomial matrix of the orthogonal Chui-Lian
  multiscaling function \cite{CL96}.
\end{itemize}

The outline of this paper is as follows.

In section \ref{sec:bauer}, we give details on Bauer's method and the
proposed modifications, and introduce the NME. In section \ref{sec:nme}
we review previous results on matrix equations. In section
\ref{sec:numerical_nme}, we discuss the details, convergence rates and
computational complexity of the numerical methods we use. In section
\ref{sec:examples}, we present some examples and verify the theoretical
convergence results from section \ref{sec:numerical_nme}.  Some
conclusions are presented in section \ref{sec:conclusion}. Some
technical results have been moved to an appendix, so as not to interrupt
the discussion in the main paper.

\section{Bauer's Method for Matrix Spectral Factorization}\label{sec:bauer}

A stable and self-correcting iteration for spectral factorization of a
para-Hermitian positive definite polynomial was proposed by F.~Bauer in
\cite{B55,B56}. Youla and Kazanjian \cite{YK78} extended this to the
positive semidefinite case.

Let $P(z)$ be a positive semidefinite para-Hermitian polynomial matrix
of degree $m$.  The equation $P(z) = H(z) H(z)^*$ is equivalent to the
fact that the infinite block banded matrix
\begin{equation}\label{eq:P}
  \calP =
  \begin{bmatrix}
    \ddots & \ddots & \ddots & \ddots & \ddots  \\
    & P_{-m} & P_{-m+1} & \cdots & P_{m-1} & P_m \\
    & & P_{-m} & P_{-m+1} & \cdots & P_{m-1} & P_m \\
    & & & \ddots & \ddots & \ddots & \ddots & \ddots
  \end{bmatrix}
\end{equation}
can be factored as $\calP = \calH \calH^*$, where $\calH$ is block
banded and lower triangular:
\begin{equation}\label{eq:H}
  \calH = 
  \begin{bmatrix}
    \ddots & \ddots & \ddots & \ddots & \\
    & H_m & \cdots & H_1 & H_0 \\
    & & H_m & \cdots & H_1 & H_0 \\
    & & & \ddots & \ddots & \ddots & \ddots & \\
  \end{bmatrix}
\end{equation}

Bauer considered the Cholesky factorization of the block banded Toeplitz
matrix $\calP_n$ obtained by selecting a $(n+1) \times (n+1)$ submatrix
of $\calP$. Thus, $\calP_n = \calH_n \calH_n^*$, where
\begin{equation}\label{eq:Pn}
  \calP_n =
  \begin{bmatrix}
    P_0    & P_1   & \cdots  & P_m & \\
    P_{-1} & P_0   & P_1   & \cdots   & P_m  & \\
    \vdots & & \ddots & & & \ddots &  \\
    P_{-m} &  & & \ddots & &  & P_m \\
    & \ddots & & & \ddots & & \vdots \\
    &  & \ddots & &  & \ddots & P_1 \\
    &  &  & P_{-m} & \cdots & P_{-1} & P_0
  \end{bmatrix}, \qquad
  \calH_n =
  \begin{bmatrix}
    H_0^{(0)}  & \\[3pt]
    H_1^{(1)} & H_0^{(1)} & \\
    \vdots & & \ddots &  & &  & \\
    H_m^{(m)} &  & & H_0^{(m)} & \\
    & \ddots & &  & \ddots &  \\
    &  & H_m^{(n)} & \cdots & H_1^{(n)}  &H_0^{(n)}
  \end{bmatrix}
\end{equation}

\subsection{Fast Bauer's Method (FBM)}
We reduce the general case to the case $m=1$, which means
$P(z) = P_{-1} z^{-1} + P_0 + P_1 z$, $H(z) = H_0 + H_1 z^{-1}$.  This
can be achieved by replacing the original coefficients with block
matrices. For $m=3$, for example, we set
\begin{equation*}
  \hat P_{-1} = \hat P_1^* =
  \begin{bmatrix}
    P_3^* & P_2^* & P_1^* \\[3pt]
    0 & P_3^* & P_2^* \\[3pt]
    0 & 0 & P_3^*
  \end{bmatrix}, \qquad
  \hat P_0 =
  \begin{bmatrix}
    P_0 & P_1 & P_2 \\[2pt]
    P_1^* & P_0 & P_1 & \\[2pt]
    P_2^* & P_1^* & P_0 
  \end{bmatrix}, \qquad
  \hat P_1 =
  \begin{bmatrix}
    P_3 & 0 & 0 \\
    P_2 & P_3 & 0 \\
    P_1 & P_2 & P_3
  \end{bmatrix}.
\end{equation*}
The resulting coefficients $\hat H_0$, $\hat H_1$ have the form
\begin{equation*}
  \hat H_1 =
  \begin{bmatrix}
    H_3 & H_2 & H_1 \\
    0 & H_3 & H_2 \\
    0 & 0 & H_3
  \end{bmatrix}, \qquad
    \hat H_0 =
  \begin{bmatrix}
    H_0 & 0 & 0 \\
    H_1 & H_0 & 0 \\
    H_2 & H_1 & H_0
  \end{bmatrix}.
\end{equation*}
Once $\hat H_0$, $\hat H_1$ have been computed, we can extract the  
original $H_0$, $H_1$, $H_2$, $H_3$ as submatrices.  

Details are given in appendix A.1.

For better readability we omit the hats on the variables in the
following, except in the statement of the algorithms.

We define $X^{(k)} = H_0^{(k)} H_0^{(k)*}$, in particular
$X^{(0)}=P_0$. Going through Cholesky's algorithm for row $(k+1)$, we
find
\begin{equation}\label{eq:fpi}
  X^{(k+1)} = P_0 - P_1^* \, [X^{(k)}]^{-1} P_1.
\end{equation}
Again, details can be found in appendix A.1.

In the limit we obtain the nonlinear matrix equation (NME)
\begin{equation}\label{eq:fixed_point}
  X = P_0 - P_1^* X^{-1} P_1.
\end{equation}
If $X$ is known, $H_0$ can be found as its Cholesky factor, and then
\begin{equation}\label{eq:findH1}
  H_1 = P_1^* H_0^{-*},
\end{equation}
where
\begin{equation*}
  H_0^{-*} = \left( H_0^{-1} \right)^* = \left( H_0^* \right)^{-1}.
\end{equation*}
Everything in the following is based on investigating the solution of
the NME \eqref{eq:fixed_point}.

Algorithm 1 summarizes our numerical approach. Step 1 in this algorithm
could be one of several numerical methods. These will be discussed in
section~\ref{sec:numerical_nme}.

It is shown in \cite{YK78} that $\calP_n = \calH_n \calH_n^*$ is
strictly positive definite for every $n$, even though the infinite
matrix $\calP$ is only positive semidefinite. This implies that all
$X^{(n)}$ are positive definite. Moreover, the limit matrix $X$ is also
positive definite, which implies that its Cholesky factor $H_0$ is
nonsingular.

It should be noted that spectral factorization is not unique, and it is
not hard to construct examples where one solution has a singular
$H_0$. However, the results in \cite{YK78} show that the particular
solution that the Bauer algorithm converges to has a nonsingular $H_0$.

\begin{center}
  \begin{tabular}{lll}
    \hline \hline
    \multicolumn{3}{c}{{\bf Algorithm 1:} A Fast Bauer's Method} \\
    \hline \hline
    \multicolumn{3}{l}{{\bf Inputs:} $P_0, P_1, \ldots, P_m$ (matrix
    coefficients of $P(z)$)} \\
    \multicolumn{3}{l}{{\bf Outputs:} $H_0, H_1, \ldots, H_m$ (matrix
    coefficients of $H(z)$)} \\
    {\bf begin} \\
    & {\bf if} $m > 1$ \\
    & & Construct block matrices $\hat P_0$, $\hat P_1$ \\
    & & to reduce to the case $m=1$ \\
    & {\bf end if} \\
    & \multicolumn{2}{l}{{\bf Step 1:} Find the matrix $\hat X$ by solving
      eq.~\eqref{eq:fixed_point}} numerically \\
    & \multicolumn{2}{l}{{\bf Step 2:} Find the matrix $\hat H_0$ as the
      Cholesky factor of $\hat X$} \\
    & \multicolumn{2}{l}{{\bf Step 3:} Find the matrix $\hat H_1$ from
      $\hat H_1 = \hat P_1^* \hat H_0^{-*}$} \\
    & {\bf if} $m > 1$ \\
    & & Extract $H_0, H_1, \ldots, H_m$ from $\hat H_0$, $\hat H_1$ \\
    & {\bf end if} \\
    {\bf end} \\
    \hline \hline
  \end{tabular}
\end{center}

\subsection{Exact Bauer's Method}

A better approach, when it works, is to solve the fixed point equation
\eqref{eq:fixed_point} in closed form. We set up the matrix $X$ with
symbolic coefficients $x_{ij}$, substitute it into
\eqref{eq:fixed_point}, and obtain a series of algebraic equations. In
simple cases, the symbolic toolbox in Matlab is able to solve this.

Based on equations \eqref{eq:fixed_point}, \eqref{eq:findH1} we obtain
Algorithm 2 for obtaining the exact values of the matrix spectral
factors.

\begin{center}
  \begin{tabular}{lll}
    \hline \hline
    \multicolumn{3}{c}{{\bf Algorithm 2:} An Exact Bauer's Method} \\
    \hline \hline
    \multicolumn{3}{l}{{\bf Inputs:} $P_0, P_1, \ldots, P_m$ (matrix
    coefficients of $P(z)$)} \\
    \multicolumn{3}{l}{{\bf Outputs:} $H_0, H_1, \ldots, H_m$ (matrix
    coefficients of $H(z)$)} \\
    {\bf begin} \\
    & {\bf if} $m > 1$ \\
    & & Construct block matrices $\hat P_0$, $\hat P_1$ \\
    & & to reduce to the case $m=1$ \\
    & {\bf end if} \\
    & \multicolumn{2}{l}{Using a suitable computer algebra system (CAS)\footnotemark,}\\
    & \multicolumn{2}{l}{{\bf Step 1:} Set up a symmetric matrix $\hat
      X$ with symbolic entries $x_{ij}$} \\
    & \multicolumn{2}{l}{{\bf Step 2:} Set up and solve the nonlinear
      system of equations} \\
    & & $f(\hat X) = \hat X - \hat P_0 + \hat P_1^* \hat X^{-1} \hat P_1 = 0$ \\
    & \multicolumn{2}{l}{Assuming that the CAS is able
      to find $\hat X$,} \\
    & \multicolumn{2}{l}{{\bf Step 3:} Find the matrix $\hat H_0$ as the Cholesky factor of $\hat X$} \\
    & \multicolumn{2}{l}{{\bf Step 4:} Find the matrix $\hat H_1$ from $\hat H_1
      = \hat P_1^* \hat H_0^{-*}$} \\
    & {\bf if} $m > 1$ \\
    & & Extract $H_0, H_1, \ldots, H_m$ from $\hat H_0$, $\hat H_1$ \\
    & {\bf end if} \\
    {\bf end} \\
    \hline \hline
  \end{tabular}
\end{center}
\footnotetext{A CAS is any mathematical software with the ability to
  solve, plot and manipulate mathematical expressions in
  analytical form.} 

\begin{reorganized}

\section{The Nonlinear Matrix Equation}\label{sec:nme}

Our goal is to find a solution of the fixed point problem
\eqref{eq:fixed_point} in the MSF setting. That is, the coefficients
$P_0$ and $P_1$ come from a para-Hermitian matrix polynomial $P(z)$
which is positive semidefinite on the unit circle. We first consider
some existence and uniqueness results for the general NME
\begin{equation}\label{eq:nme}
  X = Q - A^* X^{-1} A,
\end{equation}
where $A$, $Q \in \C^{d \times d}$, $Q > 0$, and then specialize to our
case of interest.

Equation \eqref{eq:nme} can be rephrased into two related forms. Some
results in the literature are stated for one of the other forms; in this
paper we present the equivalent results for the standard NME.

One variation is the simplified form
\begin{equation}\label{eq:nmes}
  \tilde X = I - \tilde A^* \tilde X^{-1} \tilde A,
\end{equation}
where $I$ is the identity matrix. The reduction to this form is achieved
by factoring
\begin{equation*}
  Q = M^*  M
\end{equation*}
and multiplying \eqref{eq:nme} from the left by $M^{-*}$ and from the
right by $M^{-1}$. After sorting out the terms and defining
\begin{equation}\label{eq:XPtilde}
  \tilde X = M^{-*} X M^{-1}, \qquad  \tilde A = M^{-*} A M^{-1}
\end{equation}
this becomes \eqref{eq:nmes}. It is easy to verify that $\tilde X > 0$
if and only if $X > 0$.

$M$ could be a Cholesky factor of $Q$, or we could choose $M = Q^{1/2}$
as in \cite{AMANSF08}, in which case $M = M^*$.

The other variation is the modified NME
\begin{equation}\label{eq:nme+}
  X = Q + A^* X^{-1} A.
\end{equation}
As explained in \cite{AMANSF08}, the solution $X$ of \eqref{eq:nme+} can
be computed as
\begin{equation*}
  X = Y - A Q^{-1} A^*
\end{equation*}
from the solution $Y$ of the standard NME
\begin{equation*}
  Y = R - B^* Y^{-1} B,
\end{equation*}
where $B = A Q^{-1} A$, $R = Q + A^* Q^{-1} A + A Q^{-1} A^*$.

Various kinds of nonlinear matrix equations have received much attention
in the literature \cite{ZBLDG11, LJHXZ09, ENJC93, DFGT05, JZWM09,
  CJCG10, DXLCLA11, LAYGDX10, ZXZJ96}.  Equations \eqref{eq:nme} and
\eqref{eq:nme+} were studied in \cite{ERRA93, ENJC93}, in control theory
\cite{AWNMD90, FALBC96} and in dynamic programming and statistics
\cite{ZHX96}.  The simplified NME \eqref{eq:nmes} plays a key role in
various problems of synthesis of control systems, and can also be used
to parameterize the set of stabilizing controllers \cite{AFALV08}.

If a positive definite solution exists, the NME has maximal and minimal
solutions $X_+$, $X_-$ with
\begin{equation*}
  0 < X_- \le X \le X_+
\end{equation*}
for any positive definite solution $X$.

The maximal solution $X_+$ is the unique stabilizing or weakly
stabilizing solution. A solution $X$ of NME is called {\em stabilizing}
if all the eigenvalues of the pencil $\lambda X + A$ (or equivalently
the eigenvalues of the {\em closed-loop matrix} $X^{-1} A$) are inside
the unit circle. It is {\em weakly stabilizing} if the eigenvalues are
inside or on the unit circle.

For real matrices, the NME \eqref{eq:nme} is also a special case of the
{\em Generalized Discrete Time Algebraic Riccati Equation} (GDARE)
\cite{AWFL84}
\begin{equation}\label{eq:GDARE}
  E^T X E = D^T X D - (D^T X B + A^T) (B^T X B + R)^{-1} (D^T X B + A^T)^T + C^T Q C
\end{equation}
with $R = D = 0$, $C = E = B = I$.  The matrix $X$ is a real symmetric
matrix to be found.  GDARE plays an important role in the design of
linear controllers \cite{LPRL95, MEB00}.

Some necessary and sufficient conditions for the existence of a positive
definite solution of \eqref{eq:nme} and \eqref{eq:nmes} are given in
\cite{ENJC93}.  {\em Necessary conditions} are
\begin{gather*}
  \rho(\tilde A) = \rho(A Q^{-1}) \le \frac{1}{2} \\
  \rho(\tilde A + \tilde A^*) = \rho((A + A^*) Q^{-1}) \le 1 \\
  \rho(\tilde A - \tilde A^*) = \rho((A - A^*) Q^{-1}) \le 1.
\end{gather*}
A {\em sufficient condition} is

\begin{equation*}
  \| \tilde A \|_2 = \sqrt{ \| A^* Q^{-1} A Q^{-1} \|_2} \le \frac{1}{2}.
\end{equation*}

In our MSF setting, the closed-loop matrix is $X^{-1} P_1$.  The matrix
form of the Riesz-Fej{\'e}r lemma \cite{EJL09,HHS04} guarantees the
existence of a positive definite solution. Thus, the necessary
conditions above are automatically satisfied. None of our singular
examples in section \ref{sec:examples} satisfy the sufficient condition.

The same existence result is stated in \cite[Chapter
8]{HTMLRL19}. Theorem 3.5 in \cite{ERRA93} states that this solution is
unique if and only if in addition all zeros of $|P(z)|$ are on $\T$.
\end{reorganized}

\section{Numerical Solution of the Nonlinear Matrix Equation}\label{sec:numerical_nme}

In this section, we give details on the numerical solution of
eq.~\eqref{eq:fixed_point}, which is Step 1 in Algorithm 1, with
particular attention to the singular case.

\begin{moved}

Previous algorithms
\cite{AMANSF08,FALBC96,ERRA93,ZXZJ96,ENJC93,ELSSM03,IGIVH04,CEKHLW08}
for solving NMEs require additional conditions and cannot solve spectral
factorization problems with zeroes on the unit circle very well.
\end{moved}

We consider three different approaches: fixed point iteration and
Newton's method, and numerical methods based on the eigenvalues of the
closed-loop matrix for GDARE. These methods differ in convergence speed
and accuracy, but give identical results when they work.

\begin{moved}
In addition, Matlab has a built-in nonlinear equation solver called {\tt
  fsolve}. In numerical experiments, the results of using this routine
to solve the NME were very similar to those of Newton's method. We would
expect the same to be true of other similar solvers. We will not pursue
this approach further.
\end{moved}

\subsection{Numerical Analysis Background}\label{subsec:nabg}
We consider the numerical solution of the equation
\begin{equation*}
 \vecf(\vecx) = \veczero,
\end{equation*}
where $\vecf{:} \R^d \to \R^d$. It is assumed that a solution $\vecxbar$
exists, so $\vecf(\vecxbar) = \veczero$. It is also assumed that $\vecf$
has sufficiently many continuous derivatives in a neighborhood of
$\vecxbar$.

We study stationary iterative methods of the form
\begin{align*}
  \vecx^{(0)} &= \mbox{initial guess}, \\
  \vecx^{(n+1)} &= \vecF(\vecx^{(n)}),
\end{align*}
where $\vecF$ is derived from $\vecf$ in some fashion.

We define the error in $\vecx^{(n)}$ as
\begin{equation*}
  \epsilon_\vecx^{(n)} = \| \vecxbar - \vecx^{(n)} \|
\end{equation*}
for some norm. In our numerical examples in section \ref{sec:examples}
we used the 2-norm, but other norms will give comparable results.

\begin{new}

By definition, the iteration converges if and only if
$\epsilon_\vecx^{(n)} \to 0$ as $n \to \infty$. Convergence is called
{\em linear} if
\begin{equation*}
  \epsilon_\vecx^{(n+1)} \le c \epsilon_\vecx^{(n)}
\end{equation*}
for some $c < 1$. It is called {\em quadratic} if
\begin{equation*}
  \epsilon_\vecx^{(n+1)} \le c \left[\epsilon_\vecx^{(n)}\right]^2
\end{equation*}
for some $c$. If the iteration converges, but
\begin{equation*}
  \lim_{n \to \infty}
  \frac{\epsilon_\vecx^{(n+1)}}{\epsilon_\vecx^{(n)}} =1
\end{equation*}
convergence is {\em sublinear}. Sublinear convergence can be extremely
slow; the error could behave like $1/n$, or even $1/\log n$. Typically,
the number of iterations necessary to achieve a given accuracy grows
exponentially, or even faster, with the number of digits of accuracy.

\end{new}

The residual is defined as
\begin{equation*}
  \epsilon_\vecf^{(n)} = \| \vecf(\vecxbar) - \vecf(\vecx^{(n)}) \| =  \| \vecf(\vecx^{(n)}) \|.
\end{equation*}

The beginning of the Taylor expansion of $\vecf$ around $\vecxbar$ is
\begin{equation}\label{eq:taylor}
  \vecf(\vecxbar + \Delta\vecx) \approx \vecf(\vecxbar) +
  D\,\vecf(\vecxbar) \Delta\vecx + \frac{1}{2}
  D^2\vecf(\vecxbar)(\Delta\vecx,\Delta\vecx) =
  D\,\vecf(\vecxbar) \Delta\vecx + \frac{1}{2}
  D^2\vecf(\vecxbar)(\Delta\vecx,\Delta\vecx)
\end{equation}

We call the problem {\em nonsingular} if the derivative matrix
$D\,\vecf(\vecxbar)$ is invertible, otherwise {\em singular}.

In the nonsingular case we have
\begin{align*}
  \vecf(\vecxbar + \Delta\vecx) &\approx D\,\vecf(\vecxbar) \Delta\vecx\ \\
  \|\vecf(\vecxbar + \Delta\vecx)\| &\approx \| D\,\vecf(\vecxbar) \|
  \, \|  \Delta\vecx \| = O(\|\Delta\vecx\|),
\end{align*}
so error and residual are of the same order of magnitude. We expect that
we can compute both $\vecxbar$ and $\vecf(\vecxbar + \Delta\vecx)$ to
machine accuracy, which is about 15 or 16 decimals.

In the scalar singular case we have $f'(\xbar) = 0$, therefore
\begin{align*}
  f(\xbar + \Delta x) &\approx \frac{1}{2} f''(\xbar) (\Delta x)^2, \\
  \left| f(\xbar + \Delta x) \right| &\approx \frac{1}{2} \left|
  f''(\xbar) \right| \, (\Delta x)^2 = O((\Delta x)^2).
\end{align*}
This means that the residual is approximately the square of the error.

\begin{figure}
  \centering
  \includegraphics[width=3in]{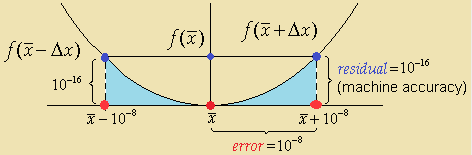}
  \caption{For $|\Delta x| = 10^{-8}$, the residual $|f(\xbar+\Delta
    x)| \approx 10^{-16}$. The function values within $10^{-8}$ of
    $\xbar$ cannot be distinguished.}
  \label{fig:1}
\end{figure}

If $|\Delta x| = 10^{-8}$, then $f(x+\Delta x) \approx 10^{-16}$, which
is machine accuracy. Since a relative change in $f(x)$ below $10^{-16}$
is not detectable, we can only compute $\xbar$ to half the available
accuracy, or about 8 decimals.  This is illustrated in Fig.~\ref{fig:1}.
For a zero of order $p$, the available accuracy in $\xbar$ reduces to
$1/p$ of the available decimals.

In the higher-dimensional singular case the situation is similar. If
$D\,\vecf(\vecxbar)$ is singular, it has eigenvectors to eigenvalue 0. A
change of order $|\Delta x|$ in the direction of one of these
eigenvectors produces a residual of order $|\Delta x|^2$. Thus, the
component of the overall error in the nullspace, and therefore the size
of the overall error, cannot be determined to more than half the
available accuracy.

\subsection{Matrix-Valued Functions}
The standard theory in subsection \ref{subsec:nabg} covers the solution
of equations of the form $\vecf(\vecx) = \veczero$, where
$\vecf{:} \R^d \to \R^d$. Eq.~\eqref{eq:fixed_point} is actually of the
form
\begin{equation*}
  f(X) =   X - P_0 + P_1^* X^{-1} P_1 = 0,
\end{equation*}
where $f$ maps $\R^{d \times d}$ to $\R^{d \times d}$.

The connection between the two cases is the $\vec$ operation \cite{T17},
which maps a matrix into a vector by stacking the columns on top of each
other. We denote the inverse operation by $\vec^{-1}$, which maps a
vector of length $d^2$ into a $d \times d$ matrix.

The matrix $X$ and matrix-valued function $f$ in $f(X) = 0$ correspond
to the vector $\vecx$ and vector-valued function $\vecf$ via
\begin{equation}\label{eq:vec}
  \vecf(\vecx) = \vec[f(\vec^{-1} \vecx)] = \veczero.
\end{equation}

The norms of $X$ and $\vecx = \vec(X)$ are comparable. In particular,
the 2-norm of $\vec(X)$ is identical to the Frobenius norm of $X$. Thus,
our concepts of error and residual are unchanged.

\subsection{Fixed Point Iteration}

Fixed point iteration is used to solve equations of the form
$\vecx = \vecg(\vecx)$. The solution $\vecxbar$ satisfies
$\vecxbar = \vecg(\vecxbar)$. This is equivalent to
$\vecf(\vecx) = \vecx - \vecg(\vecx) = \veczero$.

The iteration is
\begin{align*}
  \vecx^{(0)} &= \mbox{initial guess}, \\
  \vecx^{(n+1)} &= \vecg(\vecx^{(n)}).
\end{align*}

Bauer's method is identical to iteration \eqref{eq:fpi}, which is fixed
point iteration applied to \eqref{eq:fixed_point}, with
\begin{equation}\label{eq:g}
  g(X) = P_0 - P_1^* X^{-1} P_1.
\end{equation}
Implementing Bauer's method as a special case of fixed point iteration
does not affect the convergence speed or the result, but it is much more
efficient than the original algorithm, which forms large matrices.

\begin{new}

\subsubsection{Convergence in the Nonsingular Case}\label{subsubsec:fpi_nonsingular}
As in the Taylor expansion \eqref{eq:taylor}, we find that
\begin{equation*}
  \vecx^{(n+1)} - \vecxbar = \vecg(\vecx^{(n)}) - \vecg(\vecxbar)
  \approx D\vecg(\vecxbar) (\vecx^{(n)} - \vecxbar).
\end{equation*}
This shows that if the spectral radius of the derivative matrix
$\rho(D\vecg(\vecxbar))$ is less than 1, the iteration converges
linearly. More precisely, if $\rho(D\vecg(\vecxbar)) \le c < 1$ in a
neighborhood of $\vecxbar$, then
\begin{equation*}
  \epsilon_\vecx^{(n+1)} \le c \, \epsilon_\vecx^{(n)}.
\end{equation*}
If $\rho(D\vecg(\vecxbar)) > 1$, the iteration diverges.

\subsubsection{Convergence in the Singular Case}\label{subsubsec:fpi_singular}
In the singular case $\rho(D\vecg(\vecxbar)) = 1$ the iteration may or
may not converge. For example, in the simple scalar case $g(x) = x^2/2$
with $g'(2) = 1$ the iteration does not converge to the solution
$\xbar=2$. In the case $g(x) = 2 - (1/x)$ with $g'(1) = 1$, which is
analyzed in detail below, the iteration does converge to $\xbar=1$.

\subsubsection{Convergence for Singular MSF}\label{subsubsec:fpi_NME}
If the iteration function $\vecg$ comes from a singular MSF we have
additional results that do not apply to general $\vecg$.

Youla/Kazanjian \cite{YK78} proved that the iteration always
converges. In \cite[Algorithm 4.1]{ERRA93} it is proved that it
converges to the maximal solution $X_+$.

In the simple scalar example 2 below we can prove that the convergence
speed is sublinear, with $\epsilon^{(n)}$ of order $1/n$. This appears
to be true also for the other singular examples presented below, but
we have no proofs.
\end{new}

\subsection{Newton's Method}\label{subsec:newton}

In Newton's method, we work directly on the equation $\vecf(\vecx) =
\veczero$. The iteration is
\begin{align*}
  \vecx^{(0)} &= \mbox{initial guess}, \\
  \vecx^{(n+1)} &= \vecx^{(n)} - [D\,\vecf(\vecx^{(n)})]^{-1} \vecf(\vecx^{(n)}).
\end{align*}

\subsubsection{Implementation}
The iteration function $\vecf$ is given in \eqref{eq:vec}, but it is not
practical to calculate the derivative in this form. It is easier to
derive Newton's method for the matrix setting
\begin{equation}\label{eq:fx0}
  f(X) = X + P_1^* X^{-1} P_1 - P_0 = 0,
\end{equation}
directly.

As in eq.~\eqref{eq:ADAinv} in appendix A.2. we have
\begin{equation*}
  (X + \Delta X)^{-1} \approx X^{-1} - X^{-1} (\Delta X) X^{-1},
\end{equation*}
which leads to
\begin{equation}\label{eq:taylor2}
   \begin{split}
     f(X + \Delta X) &= (X + \Delta X) + P_1^* (X + \Delta X)^{-1} P_1 - P_0 \\
     &\approx (X + \Delta X) + P_1^* \left[ X^{-1} - X^{-1} (\Delta X) 
       X^{-1} \right] P_1 - P_0 \\
     &= f(X) + \Delta X - P_1^* X^{-1} (\Delta X) X^{-1} P_1.
   \end{split}
\end{equation}

The expression on the right is a first-order approximation to
$f(X + \Delta X)$, and Newton's method is based on setting this equal to
zero.  Thus, we compute $\Delta X$ from
\begin{equation}\label{eq:deltax}
  \Delta X - P_1^* X^{-1} (\Delta X) X^{-1} P_1 = - f(X)
\end{equation}
and use that to update $X$.

For the implementation, we use the $\vec$ operation \cite{T17}.  We
recall that for matrices $A$, $B$, $C$
\begin{equation}\label{eq:vecABC}
  \vec(ABC) = (C^T \otimes A) \, \vec(B),
\end{equation}
where $\otimes$ is the Kronecker product. Note that $C^T$ is the real
transpose, even if $C$ is a complex matrix. Eq.~\eqref{eq:deltax} is
evaluated as
\begin{equation*}
  \left[ I - (P_1^* X^{-1}) \otimes (P_1^* X^{-1}) \right] \vec(\Delta X) = - \vec(f(X)).
\end{equation*}
Using the current value $X^{(n)}$, we compute the approximate correction
$\Delta X^{(n)}$, and set $X^{(n+1)} = X^{(n)} + \Delta X^{(n)}$.

\begin{new}

\subsubsection{Convergence in the Nonsingular Case}\label{subsubsec:newton_nonsingular}
It is well known that in the nonsingular case the Newton iteration
converges quadratically to $\vecxbar$. That is,
\begin{equation*}
  \epsilon_\vecx^{(n+1)} \le c \, [\epsilon_\vecx^{(n)}]^2.
\end{equation*}

\subsubsection{Convergence in the Singular Case}\label{subsubsec:newton_singular}
In the scalar singular case convergence becomes linear. If
$f'(\xbar)=0$, but $f''(\xbar) \ne 0$, the convergence factor is
1/2. More generally, if all derivatives up to order $p$ vanish, but not
$f^{(p+1)}(\xbar)$, the convergence factor is $p/(p+1)$.

In the matrix singular case some similar results are known. We know that
the derivative matrix has an eigenvalue 0. If this eigenvalue is
nondegenerate (that is, the algebraic and geometric multiplicities are
the same), then $\R^d = M \oplus N$, where $M$ and $N$ are nullspace and
range of $D\,\vecf(\vecxbar)$, respectively. The projection of the error
onto $M$ converges linearly with factor $1/2$ (or $p/(p+1)$ in the case
of higher singularities), the projection onto $N$ converges
quadratically.  Overall, the convergence is linear. See \cite{DKK83} for
the cases $p=1$ and $p=2$, and \cite{GAOMR83} for general $p$.

If the zero eigenvalue is degenerate, the iteration may diverge or
behave chaotically \cite{GAOMR83}.

\subsubsection{Convergence for Singular MSF}\label{subsubsec:newton_NME}
It is shown in \cite{GCH99,GCH01} that in the case of a NME arising from
MSF the convergence behavior can be analyzed in terms of the spectrum of
the closed-loop matrix $X_+^{-1} P_1$, rather than the derivative
matrix. In the nonsingular case, the closed-loop matrix has all
eigenvalues inside the unit circle, and convergence is quadratic. In the
singular case there is an eigenvalue of 1. If this eigenvalue is
nondegenerate, convergence is linear with factor 1/2 (or more generally
$p/(p+1))$.

It is observed in \cite{GCH99,GCH01} and verified with multiple
numerical experiments that if the eigenvalue 1 of the closed-loop matrix
has a longest Jordan chain of length $p$, convergence of Newton's method
is linear with factor $2^{-1/p}$. However, the author of these papers
has no proof of this, nor is he aware of any proofs in the literature
\cite{GCH23}.
\end{new}

\begin{moved}

\subsection{GDARE}

Since the NME is a special case of GDARE \eqref{eq:GDARE}, algorithms
for solving GDARE could be applied. We concentrate on one which is not
based on iterative methods similar to those presented above,

Routine {\tt IDARE} in Matlab (called {\tt DARE} in earlier versions)
uses a non-iterative approach. The algorithm computes the eigenvalues of a
certain augmented matrix pencil and divides them into two groups, the
eigenvalues inside and outside the unit circle. In the singular case,
some or all eigenvalues lie on the unit circle, and are perturbed
towards the inside or outside by roundoff error. Sometimes this results
in a usable splitting, often it does not.

Convergence speed is not an issue here since there is no iteration
(other than internally to compute eigenvalues). Accuracy considerations
are the same as for fixed point iteration and Newton's method, since
these are features of the underlying problem, not a particular algorithm.

\end{moved}

\subsection{Singularity}\label{subsec:singular}

A matrix spectral factorization problem is called singular if the
determinant $\left|P(z)\right|$ has a zero on the unit circle
$\T$. Fixed point iteration is singular if the derivative matrix at the
solution has spectral radius $\rho = 1$. Newton's method is singular if
the derivative matrix at the solution is singular. In appendix A.3. we
show that these concepts are identical.

\begin{moved}

\subsection{Computational Complexity}\label{subsec:complex}

The complexity of numerical methods for Bauer's method spectral
factorizations of a matrix product filter $P(z)$ of degree $m$ with
$r \times r$ coefficients, and thus matrices $P_0$, $P_1$, $X$ of size
$mr \times mr$ is shown in Table~\ref{table:3}. Newton's method uses a
Kronecker product, so its computational complexity is the square of the
complexity of FPI or the proposed exact method. But, while Newton's
method requires only a few iterations due to its quadratic convergence,
FPI requires thousands of steps.

\begin{table}[cols=3,pos=ht]
  \caption{The complexity of numerical methods for Bauer's method
    spectral factorization of a matrix product filter $P(z)$ of degree $m$ with
$r \times r$ coefficients.}\label{table:3}
\begin{tabular}{ccc}
  \hline
  \hline
  {\bf Method} & {\bf Complexity} & {\bf Iterations} \\
  \hline
  \hline
  Newton & $O((mr)^6)$ & a few \\
  FPI & $O((mr)^3)$ & thousands \\
  \hline
  \hline
\end{tabular}
\end{table}
\end{moved}

\section{Numerical Examples}\label{sec:examples}

We illustrate our approach with several examples, taken from different
areas: control theory (example 1), a simple scalar example (example 2),
a well-known example of Lasha Ephremidze (example 3), and orthogonal
multiwavelet filter banks that have played a significant role in signal
processing (examples 4--7).

The examples are arranged in order of increasing singularity and
difficulty. Example 1 is nonsingular. Example 2 is singular, but scalar,
and can be completely analyzed. Example 3 has two double zeros on
$\T$. Examples 4--6 have quadruple zeros, and Example 7 has a tenfold
zero (and is of larger size $5 \times 5$).

The considered examples are tabulated in Table~\ref{table:1}, which
shows the singularity of polynomial matrices and zeros on the unit
circle.

\begin{table}[cols=3,pos=ht]
  \caption{The characteristics of the scalar case (example 2) and the
    matrix cases (examples 1, 3--7).}\label{table:1}
\begin{tabular}{ccc}
  \hline
  \hline
  Examples & Singularity & Zeros on the unit circle \\
  \hline
  \hline
  1 & No & none \\
  2 & Yes	& two double \\
  3 & Yes	& two double \\
  4 & Yes & quadruple \\
  5 & Yes & quadruple \\
  6 & Yes & quadruple \\
  7 & Yes & decuple \\
  \hline
  \hline
\end{tabular}
\end{table}

\subsection{The Nonsingular Case}

{\bf Example 1: An example from Ku{\v c}era \cite[p.~181]{K79}}

We consider the para-Hermitian matrix polynomial of order $m = 2$
\begin{equation*}
  P(z) =
  \begin{bmatrix}
    0 & 0 \\
    2 & 0
  \end{bmatrix} z^{-2} + 
  \begin{bmatrix}
    0 & 0 \\
    0 & -2
  \end{bmatrix} z^{-1} +
  \begin{bmatrix}
    1 & 0 \\ 
    0 & 9
  \end{bmatrix} +
  \begin{bmatrix}
    0 & 0 \\
    0 & -2
  \end{bmatrix} z +
  \begin{bmatrix}
    0 & 2 \\
    0 & 0
  \end{bmatrix} z^2.
\end{equation*}
The determinant $|P(z)| = z^{-1}(2z - 1)(2 - z)$ has zeros at 1/2 and
2, so this matrix polynomial is nonsingular. We need to form block
matrices $\hat P_0$, $\hat P_1$ of size $4 \times 4$.

This problem is too large for a symbolic solution, but Algorithm 1 with
Newton's method converges extremely rapidly and accurately to a spectral
factor whose entries can be recognized as rational. We find
\begin{equation}\label{eq:ex1sol}
  \hat X = \frac{1}{17}
  \begin{bmatrix}
    8 & 2 & -2 & 0 \\
    2 & 145 & 8 & -34 \\
    -2 & 8 & 9 & 0 \\
    0 & -34 & 0 & 153
  \end{bmatrix}, \quad
  H(z) = \frac{1}{\sqrt{34}}
  \begin{bmatrix}
    4 & 0 \\
    1 & 17
  \end{bmatrix} + \frac{1}{\sqrt{34}}
  \begin{bmatrix}
    -1 & 1 \\
    0 & -4
  \end{bmatrix} z^{-1} + \frac{1}{\sqrt{34}}
  \begin{bmatrix}
    0 & 4 \\
    0 & 0
  \end{bmatrix} z^{-2}.
\end{equation}

\subsection{The Scalar Singular Case}

{\bf Example 2: The (non-normalized) scalar product filter of the Haar scaling function}

We consider the scalar polynomial with $r=1$, $m=1$
\begin{equation*}
  P(z) = z^{-1} + 2 + z = (1+z)(1+z^{-1}) = \frac{(z+1)^2}{z}.
\end{equation*}
The determinant is $P(z)$ itself, which has a double zero at $z=(-1)$,
so this is a singular example.

The solution is
\begin{equation*}
  X = 1, \quad H(z) = 1 + z^{-1}.
\end{equation*}

This example is included because its convergence speed can be completely
analyzed.

\subsection{The Matrix Singular Case}

The following examples are sorted by increasing singularity, from a
double zero to a tenth order zero on the unit circle.

\noindent {\bf Example 3: An example from Ephremidze et al.~\cite{ESS18}}

We consider the singular polynomial matrix of order $m=1$
\begin{equation*}
  P(z) =
  \begin{bmatrix}
    2 & 11 \\
    7 & 38
  \end{bmatrix} z^{-1} +
  \begin{bmatrix}
    6 & 22 \\
    22 & 84
  \end{bmatrix} +
  \begin{bmatrix}
    2 & 7 \\
    11 & 38
  \end{bmatrix} z.
\end{equation*}
The determinant $|P(z)| = -(z+1)^2 (z-1)^2 / z^2$ has two double zeros
($z = \pm 1$) on the unit circle.

The solution can be computed in closed form. It is
\begin{equation*}
  X = 
  \begin{bmatrix}
    1 & 5 \\
    5 & 26
  \end{bmatrix}, \quad
  H(z) =
  \begin{bmatrix}
    1 & 0 \\
    5 & 1
  \end{bmatrix} +
  \begin{bmatrix}
    2 & 1 \\
    7 & 3
  \end{bmatrix} z^{-1},
\end{equation*}
which is the same as the result in \cite{ESS18}.

\noindent {\bf Example 4: An integer multiwavelet filter bank from Cheung et al.\ \cite{CCP99}}

We consider the para-Hermitian polynomial matrix of order $m=1$
\begin{equation*}
  P(z) = \frac{1}{4}
  \begin{bmatrix}
    2 & -\sqrt{2} \\
    \sqrt{2} & 0
  \end{bmatrix} z^{-1} + 
  \begin{bmatrix}
    1 & 0 \\
    0 & 1
  \end{bmatrix} + \frac{1}{4}
  \begin{bmatrix}
    2 & \sqrt{2} \\
    -\sqrt{2} & 0
  \end{bmatrix} z
\end{equation*}
The determinant $|P(z)| = (1 + z)^4/(8z^{-2})$ has a quadruple zero at
$z = -1$.

The solution can be computed in closed form. It is
\begin{equation*}
  X = \frac{1}{4}
  \begin{bmatrix}
    2 & - \sqrt{2} \\
    -\sqrt{2} & 2
  \end{bmatrix}, \quad
  H_0 = \frac{1}{2}
  \begin{bmatrix}
    \sqrt{2} & 0 \\
    -1 & 1
  \end{bmatrix}, \quad H_1 = \frac{1}{2}
  \begin{bmatrix}
    \sqrt{2} & 0 \\
    1 & 1
  \end{bmatrix}.
\end{equation*}

The spectral factor $H(z)$ represents an orthogonal multiscaling
function, which can be completed with the multiwavelet function
\begin{equation*}
  G_0 = \frac{1}{2}
  \begin{bmatrix}
    0 & \sqrt{2} \\
    1 & 1
  \end{bmatrix}, \quad G_1 = \frac{1}{2}
  \begin{bmatrix}
    0 & -\sqrt{2} \\
    -1 & 1
  \end{bmatrix}.
\end{equation*}

After left multiplying with the matrix $C = \diag(\sqrt{2},1)$ we obtain
the well-known multiwavelet with dyadic coefficients
\begin{gather*}
  H_0 =
  \begin{bmatrix}
    1 & 0 \\
    -\frac{1}{2} & \frac{1}{2}
  \end{bmatrix}, \quad H_1 = 
  \begin{bmatrix}
    1 & 0 \\
    \frac{1}{2} & \frac{1}{2}
  \end{bmatrix}, \quad G_0 =
  \begin{bmatrix}
    0 & 1 \\
    \frac{1}{2} & \frac{1}{2}
  \end{bmatrix}, \quad G_1 = 
  \begin{bmatrix}
    0 & -1 \\
    -\frac{1}{2} & \frac{1}{2}
  \end{bmatrix},
\end{gather*}
as described in \cite{CCP99}. The multiwavelet decomposition maps dyadic
numbers to dyadic numbers, which guarantees a perfect reconstruction in
the absence of quantizers. It leads to lossless processing, and its
lifting scheme can be easily implemented using hardware or software.

\noindent {\bf Example 5: A new supercompact multiwavelet filter}

From the Chui-Lian multiscaling function \cite{CL96} we obtain the
singular matrix product filter
\begin{equation*}
  P(z) = \frac{1}{8}
  \begin{bmatrix}
    4 & 1 + \sqrt{7} \\ 
    -(1 + \sqrt{7}) & -\sqrt{7}
  \end{bmatrix} z^{-1}+
  \begin{bmatrix}
    1 & 0 \\
    0 & 1
  \end{bmatrix} +
  \frac{1}{8}
  \begin{bmatrix}
    4 & -(1 + \sqrt{7}) \\
    1 + \sqrt{7} & -\sqrt{7}
  \end{bmatrix} z.
\end{equation*}
The determinant $|P(z)| = (4-\sqrt{7})(1+ z)^4/(32 z^2 )$ has a
quadruple zero at $z = -1$. 

The solution can be computed in closed form. It is
\begin{equation*}
  X = \frac{1}{8}
  \begin{bmatrix}
    4 & \sqrt{7}+1 \\
    \sqrt{7}+1 & 4
  \end{bmatrix}, \quad
  H_0 = \frac{\sqrt{2}}{8}
  \begin{bmatrix}
    4 & 0 \\
    \sqrt{7} + 1 & \sqrt{7} - 1
  \end{bmatrix}, \quad H_1 = \frac{\sqrt{2}}{8}
  \begin{bmatrix}
    4 & 0 \\
    -(\sqrt{7}+1) & \sqrt{7} - 1
  \end{bmatrix}.
\end{equation*}
As in example 4, the spectral factor $H(z)$ represents an orthogonal
multiscaling function, which can be completed with the multiwavelet
function
\begin{equation*}
  G_0 = \frac{\sqrt{2}}{8}
  \begin{bmatrix}
    0 & 4 \\
    \sqrt{7} - 1 & -\sqrt{7} - 1
  \end{bmatrix}, \quad G_1 = \frac{\sqrt{2}}{8}
  \begin{bmatrix}
    0 & -4 \\
    -\sqrt{7}+1 & -\sqrt{7} - 1
  \end{bmatrix}.
\end{equation*}

The determinant of $H(z)$ includes the factor $(1 + z^{-1})^2$, which
implies regularity. This multiwavelet filter is supercompact (i.\ e.\ it
has support $[0,1]$). Despite its short length, this multifilter has
better Sobolev smoothness \cite{J98} ($S_{MSF} = 1.28$) than the
Chui-Lian multiwavelet filter ($S_{CL}=1.06$).

\noindent {\bf Example 6: Legendre multiscaling function of order 2}

We consider the singular polynomial matrix of order $m=1$ obtained in
\cite{KCK20}
\begin{equation*}
  P(z) = \frac{1}{4}
  \begin{bmatrix}
    2 & \sqrt{3} \\
    -\sqrt{3} & -1
  \end{bmatrix} z^{-1} +
  \begin{bmatrix}
    1 & 0 \\
    0 & 1
  \end{bmatrix} + \frac{1}{4}
  \begin{bmatrix}
    2 & -\sqrt{3} \\
    \sqrt{3} & -1
  \end{bmatrix} z.
\end{equation*}
The determinant $|P(z)| = (z+1)^4/(16 z^2)$ has a quadruple zero at
$z = -1$.

The solution can be computed in closed form. It is
\begin{equation*}
  X = \frac{1}{4}
  \begin{bmatrix}
    2 & \sqrt{3} \\
    \sqrt{3} & 2
  \end{bmatrix}, \quad
  H_0 = \frac{\sqrt{2}}{4}
  \begin{bmatrix}
    2 & 0 \\
    \sqrt{3} & 1
  \end{bmatrix}, \quad
  H_1 = \frac{\sqrt{2}}{4}
  \begin{bmatrix}
    2 & 0 \\
    -\sqrt{3} & 1
  \end{bmatrix}.
\end{equation*}
This again corresponds to an orthogonal symmetric/antisymmetric
multiscaling function.

\noindent {\bf Example 7: Legendre multiscaling function of order 5 \cite{JGZ15}}

This is again a para-Hermitian polynomial matrix of order $m=1$ which
forms a supercompact multiwavelet filter bank, but with matrices of size
$5 \times 5$:
\begin{equation*}
  P_1 = \frac{1}{256}
  \begin{bmatrix}
    128 & 64 \sqrt{3} & 0 & -16 \sqrt{7} & 0 \\
    -64 \sqrt{3} & -64 & 16 \sqrt{15} & 16 \sqrt{21} & -8 \sqrt{3} \\
    0 & -16 \sqrt{15} & -112 & -8 \sqrt{35} & 24 \sqrt{5} \\
    16 \sqrt{7} & 16 \sqrt{21} & 8 \sqrt{35} & -40 & -39 \sqrt{7} \\
    0 & 8 \sqrt{3} & 24 \sqrt{5} & 39 \sqrt{7} & 53
  \end{bmatrix}, \quad P_0 =
  \begin{bmatrix}
    1 & 0 & 0 & 0 & 0 \\
    0 & 1 & 0 & 0 & 0 \\
    0 & 0 & 1 & 0 & 0 \\
    0 & 0 & 0 & 1 & 0 \\
    0 & 0 & 0 & 0 & 1
  \end{bmatrix}.
\end{equation*}

The determinant $|P(z)| = (z+1)^{10}/(2^{25}\,z^5)$ has a tenfold zero
at $z=-1$.

The original system is too large to be solved in closed form. However,
the numerical solution reveals certain relationships between the
coefficients in $X$. For example, it becomes apparent where zeroes are
located, or which coefficients are multiples of others. Using this
information to reduce the number of variables involved leads to a
problem that can be solved in closed form.

We obtain the spectral factor
\begin{gather*}
  X = \frac{1}{256}
  \begin{bmatrix}
    128 & -64\sqrt{3} & 0 & 16\sqrt{7} & 0 \\
    -64\sqrt{3} & 128 & -16\sqrt{15} & 0 & 8\sqrt{3} \\
    0 & -16\sqrt{15} & 128 & -16\sqrt{35} & 0 \\
    16\sqrt{7} & 0 & -16\sqrt{35} & 128 & -21\sqrt{7} \\
    0 & 8\sqrt{3} & 0 & -21\sqrt{7} & 128
  \end{bmatrix}, \\[10pt]
  H_0 = \frac{\sqrt{2}}{32}
  \begin{bmatrix}
    16 & 0 & 0 & 0 & 0 \\
    -8 \sqrt{3} & 8 & 0 & 0 & 0 \\
    0 & -4 \sqrt{15} & 4 & 0 & 0 \\
    2 \sqrt{7} & 2 \sqrt{21} & -2 \sqrt{35} & 2 & 0 \\
    0 & 2 \sqrt{3} & 6 \sqrt{5} & 3 \sqrt{7} & 1
  \end{bmatrix}, \quad
  H_1 = \frac{\sqrt{2}}{32}
  \begin{bmatrix}
    16 & 0 & 0 & 0 & 0 \\
    8 \sqrt{3} & 8 & 0 & 0 & 0 \\
    0 & 4 \sqrt{15} & 4 & 0 & 0 \\
    - 2 \sqrt{7} & 2 \sqrt{21} & 2 \sqrt{35} & 2 & 0 \\
    0 & - 2 \sqrt{3} & 6 \sqrt{5} & -3 \sqrt{7} & 1
  \end{bmatrix}.
\end{gather*}
Again, these correspond to a symmetric/antisymmetric multiscaling function.

\subsection{Convergence Rate Analysis for Fixed Point Iteration}

The error in fixed point iteration can be measured in two different
ways: By the error in the spectral factor
\begin{equation*}
  \epsilon_H^{(n)} = \| H_0 - H_0^{(n)} \|
\end{equation*}
(which we can only calculate if the true coefficient $H_0$ is known to
sufficient accuracy), or by the residual
\begin{equation*}
  \epsilon_P^{(n)} = \| P_0 - H_0^{(n)} H_0^{(n)*} - H_1^{(n)} H_1^{(n)*} \|.
\end{equation*}
It is shown in appendix A.2.~that it is not necessary to consider the
error in $H_1$ or the residual in $P_1$.

In this section we determine the convergence speed and accuracy of fixed
point iteration for the seven examples and compare them to the
theoretical considerations from section \ref{sec:numerical_nme}.

\subsubsection{The Nonsingular Case}

The convergence behavior in example 1 is exactly as
expected. Fixed-point iteration converges linearly, that is,
$\epsilon^{(n+1)} \approx c \, \epsilon^{(n)}$, so
$\epsilon^{(n)} \approx c^n \, \epsilon^{(0)}$ for both $\epsilon_H$ and
$\epsilon_P$.

\begin{figure}
  \centering
  \includegraphics[width=6in]{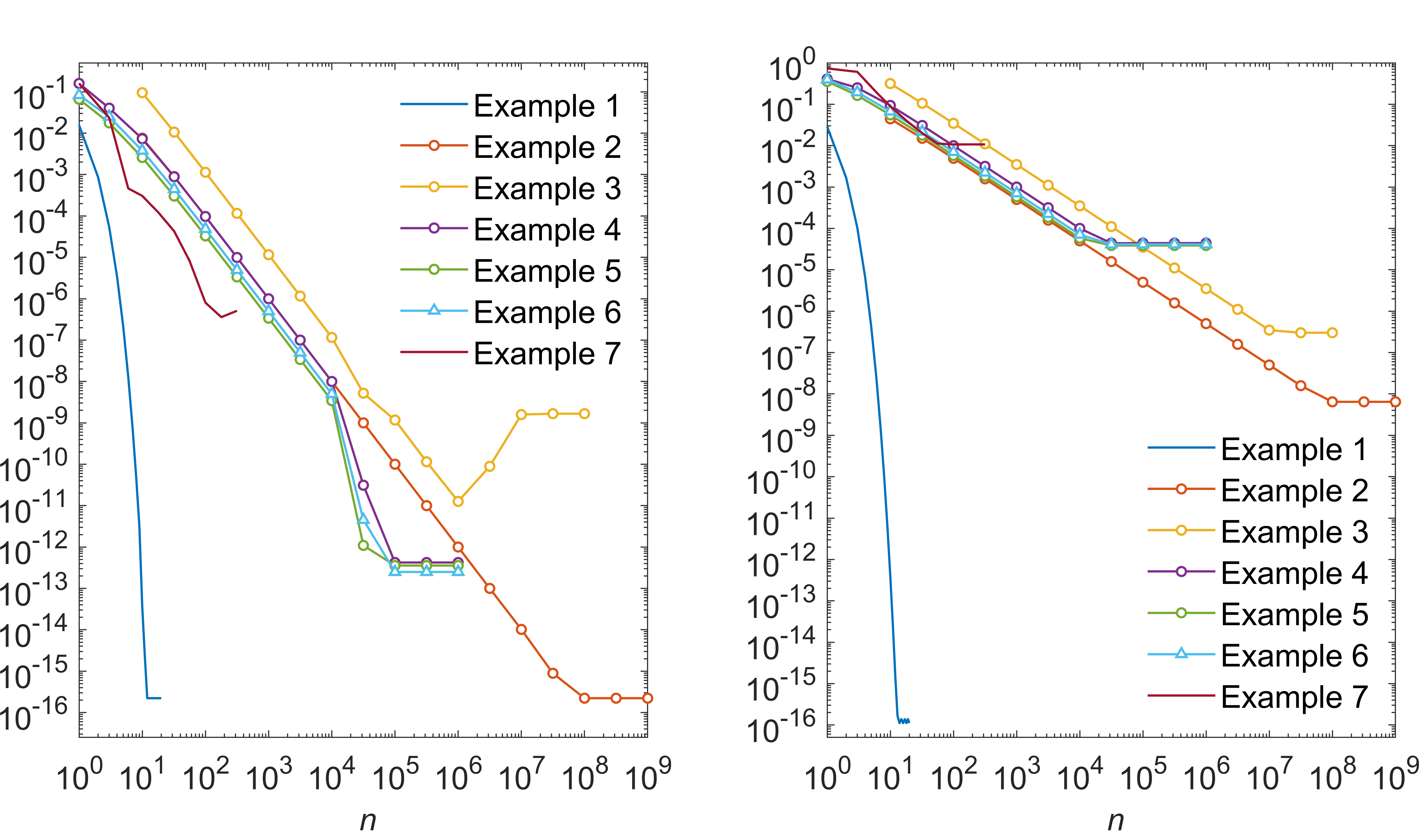} \\
  (a) \hspace{2.6in} (b)
  \caption{Log-log plots of numerical errors obtained by {\em Fixed
      Point Iteration} (FPI) applied to Fast Bauer’s method (FBM) for
    scalar and matrix spectral factorization; (a) the residuals
    $\epsilon_P$ of the product filter; (b) the errors $\epsilon_H$ of
    the spectral factor.}
  \label{fig:fpi}
\end{figure}

Figure \ref{fig:fpi} illustrates that $\epsilon_P = O(\epsilon_H)$ and
that convergence stops once we reach machine accuracy.

In comparison, in \cite{K79} this method with a stopping criterion of
the error $\epsilon = 10^{-4}$ leads to a spectral factor which is equal
to the exact solution \eqref{eq:ex1sol} to four decimals
\begin{equation*}
  H(z) =
  \begin{bmatrix}
    0.6860 & 0 \\
    0.1715 & 2.9155
  \end{bmatrix} + 
  \begin{bmatrix}
    -0.1715 & 0.1715 \\
    0 & -0.6860 
  \end{bmatrix} z^{-1} + 
  \begin{bmatrix}
    0 & 0.6860 \\
    0 & 0
  \end{bmatrix} z^{-2}
\end{equation*}
which is the true solution to four decimals.

Algorithm 1 in \cite{MTJ09}, using 30 steps of size 0.23, leads to a
spectral factor with big errors (shown in bold):
\begin{equation*}
  H(z) =
  \begin{bmatrix}
    0.6860 & \mathbf{0.0326} \\
    \mathbf{0.0326} & 2.9155 
  \end{bmatrix} +
  \begin{bmatrix}
    -0.1715 & 0.1715 \\
    \mathbf{0.0326} & -0.6860 
  \end{bmatrix} z^{-1} + 
  \begin{bmatrix}
    0 & 0.6860 \\
    0 & 0
  \end{bmatrix} z^{-2}.
\end{equation*}

\subsubsection{The Scalar Singular Case}

The fixed point iteration in example 2 is
\begin{equation*}
  x^{(n+1)} = 2 - \frac{1}{x^{(n)}},
\end{equation*}
which should converge to $\xbar = 1$. The error
$\epsilon_H^{(n)} = x^{(n)} - 1$ satisfies
\begin{equation*}
  \epsilon_H^{(n+1)} = \frac{\epsilon_H^{(n)}}{\epsilon_H^{(n)}+1}.
\end{equation*}
This is one-sided convergence: if $\epsilon_H^{(k)} > 0$ for some $k$,
then $\epsilon_H^{(n)} > 0$ for all $n>k$, and the errors are monotone
decreasing. In our case $x^{(0)} = 2$, $e^{(0)} = 1$, so this is
satisfied.

Some algebra shows
\begin{equation*}
  \epsilon_H^{(1)} = \frac{\epsilon_H^{(0)}}{\epsilon_H^{(0)} + 1}, \qquad 
  \epsilon_H^{(2)} = \frac{\epsilon_H^{(0)}}{2\epsilon_H^{(0)} + 1}, \quad
  \cdots \quad \epsilon_H^{(n)} = \frac{\epsilon_H^{(0)}}{n\epsilon_H^{(0)} + 1},
\end{equation*}
so
\begin{equation*}
  \epsilon_H^{(n)} \approx \frac{1}{n}, \qquad \ln \epsilon_H^{(n)} \approx - \ln n.
\end{equation*}
This is sublinear convergence. A plot of $\ln \epsilon_H$ versus $\ln n$
should be a line with slope $(-1)$.  This is illustrated in
fig.~\ref{fig:fpi}(b).

We can also calculate that
\begin{equation*}
  H_1^{(n)} = \frac{P_1}{H_0^{(n)}} = \frac{1}{1 + \epsilon_H^{(n)}}
  \approx 1 -  \epsilon_H^{(n)},
\end{equation*}
so the residual in $P_0$ is
\begin{align*}
  \epsilon_P^{(n)} &= [H_0^{(n)}]^2 + [H_1^{(n)}]^2 - P_0 \\
  &\approx (1 + \epsilon_H^{(n)})^2 +  (1 - \epsilon_H^{(n)})^2 - 2 =
  O([\epsilon_H^{(n)}]^2).
\end{align*}

A plot of $\ln \epsilon_P$ versus $\ln n$ should be a line with slope
$(-2)$, which is illustrated in fig.~\ref{fig:fpi}(a).

\subsubsection{The Matrix Singular Case}

In example 3, the determinant has double zeros on $\T$. The achievable
accuracy for the spectral factors, for both fixed point iteration and
Newton's method, is approximately $10^{-8}$, which is half the number of
available digits. Examples 4--6 have quadruple zeros. The convergence
stalls near accuracy $10^{-4}$. Example 7 has a tenfold zero. The
convergence stalls near accuracy $10^{-1}$ or $10^{-2}$. These are all
the same results we expect in the scalar case.

For fixed point iteration, we observe the same kind of sublinear
convergence as in the scalar example:
\begin{equation*}
  \epsilon_H^{(n)} = O(\frac{1}{n}), \quad 
  \epsilon_P^{(n)} = O([\epsilon_H^{(n)}]^2).
\end{equation*}

In some cases we see strange effects if we try to push beyond the point
of highest accuracy. Convergence becomes erratic, and it appears that
$\epsilon_P^{(n)} = O([\epsilon_H^{(n)}]^3)$ or even
$\epsilon_P^{(n)} = O([\epsilon_H^{(n)}]^4)$. These numerical
instabilities are due to ill-conditioning of the underlying problem and
will occur for any algorithm.

\subsection{Convergence Rate Analysis for Newtons' Method}

\subsubsection{The Nonsingular Case}

\begin{figure}
   \centering
  \includegraphics[width=6in]{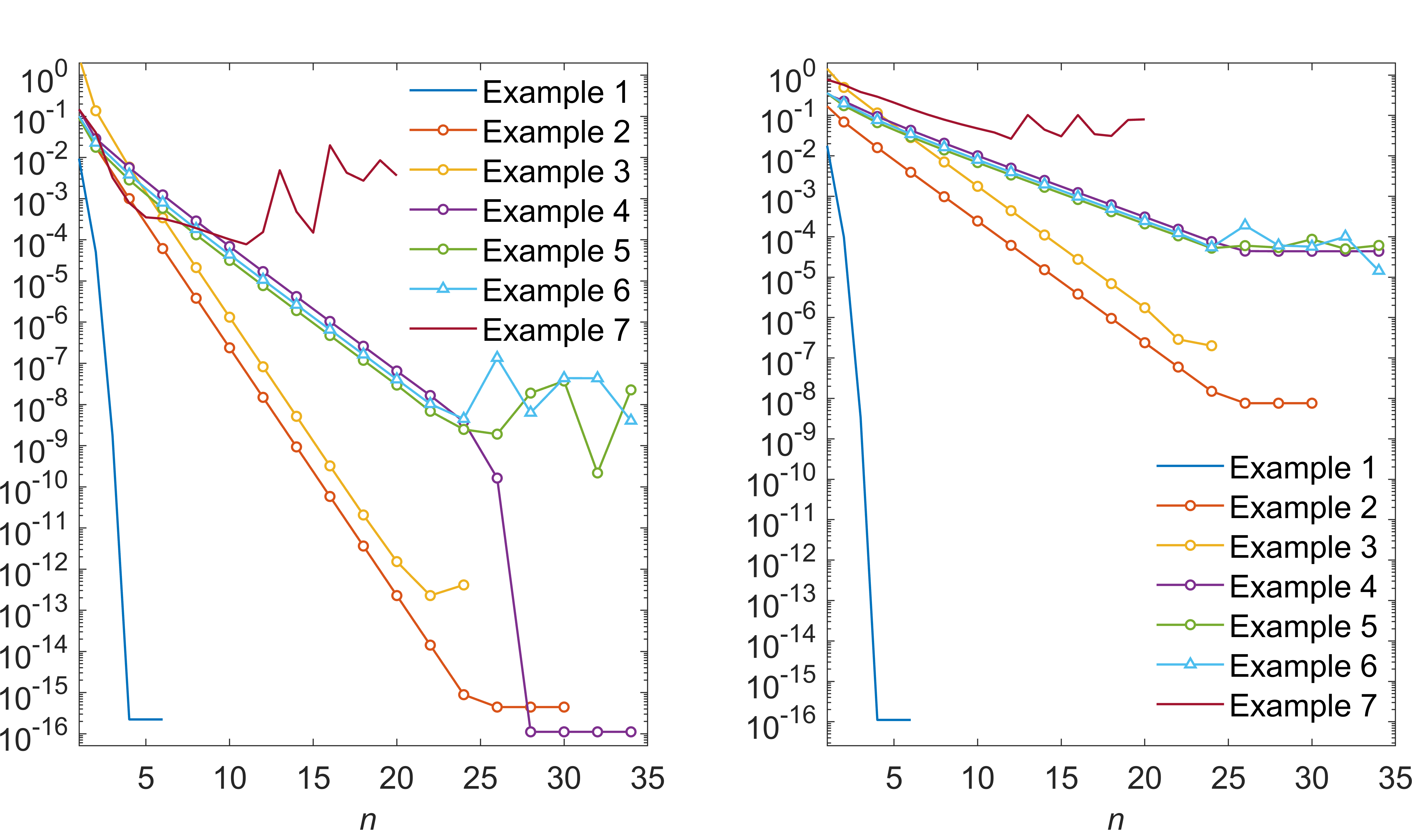} \\
  (a) \hspace{3in} (b)
  \caption{Semi-log plots of numerical errors obtained by {\em Newton's
      Method} applied to Fast Bauer’s method (FBM) for scalar and matrix
    spectral factorization; (a) the residuals $\epsilon_P$ of the
    product filter; (b) the errors $\epsilon_H$ of the spectral factor.}
   \label{fig:newton}
 \end{figure}

Fast Bauer’s Method using Newton's method in example 1, which is quadratically
convergent, achieves full precision with errors
$\epsilon_H = 1.11 \cdot 10^{-16}$ and
$\epsilon_P = 2.22 \cdot 10^{-16}$ with only 4 iterations (see
fig.~\ref{fig:newton}(a) and (b)).

Figure \ref{fig:newton} illustrates that $\epsilon_P = O(\epsilon_H)$
and that convergence stops once we reach machine accuracy.

\subsubsection{The Scalar Singular Case}

Newton's method in the scalar case is known to converge linearly with
factor $p/(p+1)$ if $f(x)$ has a $p$-fold zero at $\xbar$.  Moreover,
the residual is the $p^{th}$ power of the error, so we can only compute
the solution to $1/p$ times the available digits.

This is easy to verify directly in example 2. Equation \eqref{eq:fx0} in
this case becomes
\begin{equation*}
  f(x) = x - p_0 + p_1^* x^{-1} p_1 = x - 2 + 1 \cdot x^{-1} \cdot 1 = 0,
\end{equation*}
so 
\begin{equation*}
  f(x) = x - 2 + \frac{1}{x}.
\end{equation*}
This leads to
\begin{equation*}
  x^{(n+1)} = 2 - \frac{2}{x^{(n)} + 1}, \qquad \epsilon_H^{(n+1)} = \frac{\epsilon_H^{(n)}}{\epsilon_H^{(n)}+2}.
\end{equation*}
Some algebra gives
\begin{equation*}
  \epsilon_H^{(1)} = \frac{\epsilon_H^{(0)}}{\epsilon_H^{(0)}+2}, \qquad 
  \epsilon_H^{(2)} = \frac{\epsilon_H^{(0)}}{3 \epsilon_H^{(0)}+4}, \qquad
  \epsilon_H^{(3)} = \frac{\epsilon_H^{(0)}}{7 \epsilon_H^{(0)}+8}, \quad
  \cdots \quad
  \epsilon_H^{(n)} = \frac{\epsilon_H^{(0)}}{(2^n-1) \epsilon_H^{(0)}+2^n},
\end{equation*}
so
\begin{equation*}
  \epsilon_H^{(n)} \approx 2^{-n} \epsilon_H^{(0)}.
\end{equation*}
As in the case of fixed point iteration, we can verify directly that
$\epsilon_P^{(n)} = O([\epsilon_H^{(n)}]^2)$.

If we plot $\ln(\epsilon_H)$ versus $n$, this should be approximately a
straight line with a slope of $(- \ln 2)$. A plot of $\ln(\epsilon_P)$
should be a straight line with slope $(-2\ln 2)$.  This is illustrated
in fig.~\ref{fig:newton}(a) and (b).  Notice how $\epsilon_H$ stalls
around $10^{-8}$.

\subsubsection{The Matrix Singular Case}
\begin{new}

For Newton's method, we observe
\begin{equation*}
  \epsilon_H^{(n)} \approx c^n \epsilon_0^{(n)}, \quad
  \epsilon_P^{(n)} = O([\epsilon_H^{(n)}]^2).
\end{equation*}
The behavior corresponds with the convergence results quoted in
subsection \ref{subsubsec:newton_singular}.

In example 3, $c = 1/2$, since the eigenvalue 1 of the closed-loop
matrix is nondegenerate. In examples 4--6, $c=1/\sqrt{2} \approx 0.7$,
since the eigenvalue 1 has a Jordan chain of length 2. In example 7, $c
= 1/\sqrt[5]{2} \approx 0.87$, since the longest Jordan chain has length
5. 

Convergence becomes even more erratic than for fixed point iteration if
you try to push beyond the highest achievable accuracy.
\end{new}

Another numerical problem arises in example 7.  For $k > 64$ the matrix
$P^{(k)}$ is not positive definite any more, due to roundoff error; the
algorithm cannot be continued past this point. After $k = 40$, the
residual error $\epsilon_P$ decreases, while the error $\epsilon_H$
increases and ends up near $10^{-2}$. This shows the strong influence of
a high order zero on the precision of the spectral factors.

\begin{reorganized}

\subsection{Results of Standard GDARE Routines}

We consider the NME as a special case of GDARE \eqref{eq:GDARE} and use
existing numerical routines for this.

One such routine is {\tt DiscreteRiccatiSolve} in Mathematica.
Unfortunately, it requires $R$ in \eqref{eq:GDARE} to be nonsingular,
while $R=0$ in our case. This routine cannot be applied.

Maple has a routine {\tt DARE} which can solve all of the examples
except example 3. Results are shown in table~\ref{table:maple}.
According to the obtained results, the errors $\epsilon_X$, $\epsilon_H$
and $\epsilon_P$ for examples 1--3 are satisfactory, while examples 4--7
have very low precision.

\begin{table}[cols=5,pos=ht]
  \caption{The errors $\epsilon_X$, $\epsilon_H$ and $\epsilon_P$ for
    examples 1--7, using the DARE routine in Maple 17.}\label{table:maple}
\begin{tabular}{c||ccc}
  Example & $\epsilon_X$ & $\epsilon_H$ & $\epsilon_P$ \\
  \hline
  \hline
  1 & 2.35e-11 & 7.01e-11 & 3.29e-9 \\
   \cline{2-4}
  \multirow{4}{*}{2} & \multicolumn{3}{c}{$p_0 = 1$, $p_1 = 1/2$} \\
               & 1.22e-8 & 6.10e-9 & 0 \\
               & \multicolumn{3}{c}{$p_0 = 2$, $p_1 = 1$} \\
               & 6.10e-9 & 4.50e-9 & 5.55e-17 \\
   \cline{2-4}
               3 & \multicolumn{3}{c}{Incorrect solution $X$} \\
               4 & 5.95e-4 & 5.95e-4 & 6.59e-10 \\
               5 & 1.69e-4 & 2.90e-4 & 6.59e-10 \\
               6 & 2.20e-4 & 5.95e-4 & 2.64e-10 \\
               7 & 2.61e-3 & 2.61e-3 & 2.61e-3 \\
  \hline
  \hline
\end{tabular}
\end{table}

Routine {\tt IDARE} in Matlab (called {\tt DARE} in earlier versions)
can solve nonsingular GDARE problems, but its performance on singular
equations is rather spotty. Numerical results for Matlab versions
R2012a, R2018a, R2020a, R2020b, R2021a and R2021a Update 4 are shown in
Table~\ref{table:matlab}. The results are different for different versions.

GDARE \eqref{eq:GDARE} is associated with the matrix pencil
\begin{equation*}
  M - zN =
  \begin{bmatrix}
    R & 0 & B \\
    C^T Q C & E^T & -A^T \\
    A & 0 & R
  \end{bmatrix} - z
  \begin{bmatrix}
    E & 0 & 0 \\
    0  & A^T  & 0 \\
    0 & -B^T & 0
  \end{bmatrix}.
\end{equation*}
If we impose $E = B = C = I$, $R = D = 0$, $Q = P_0 = I$ and $A = P_1$,
the matrix pencil is
\begin{equation*}
  M - zN =
  \begin{bmatrix}
    0 & 0 & I \\
    I & I & -P_1^T \\
    P_1 & 0 & 0
  \end{bmatrix} - z
  \begin{bmatrix}
    I & 0 & 0 \\
    0  & 0  & 0 \\
    0 & -I & 0
  \end{bmatrix}.
\end{equation*}
The eigenvalues of this matrix pencil, or equivalently the eigenvalues
of the closed loop matrix, are those chosen by the {\tt IDARE}
algorithm. In the singular case, some of them lie on the unit circle and
roundoff error can perturb them to the outside.

This is already obvious for the very simple Example 2. None of the IDARE
versions were able to solve it, but some versions work if the example is
simply scaled by a factor 1/2.

Fig.~\ref{fig:are1} shows these eigenvalues for our examples. In example
1 they lie inside the unit circle; in example 7 they lie either in or on
the unit circle. For examples 2--6 all of them lie on the unit
circle.

\begin{figure}
  \centering
  \includegraphics[width=3in]{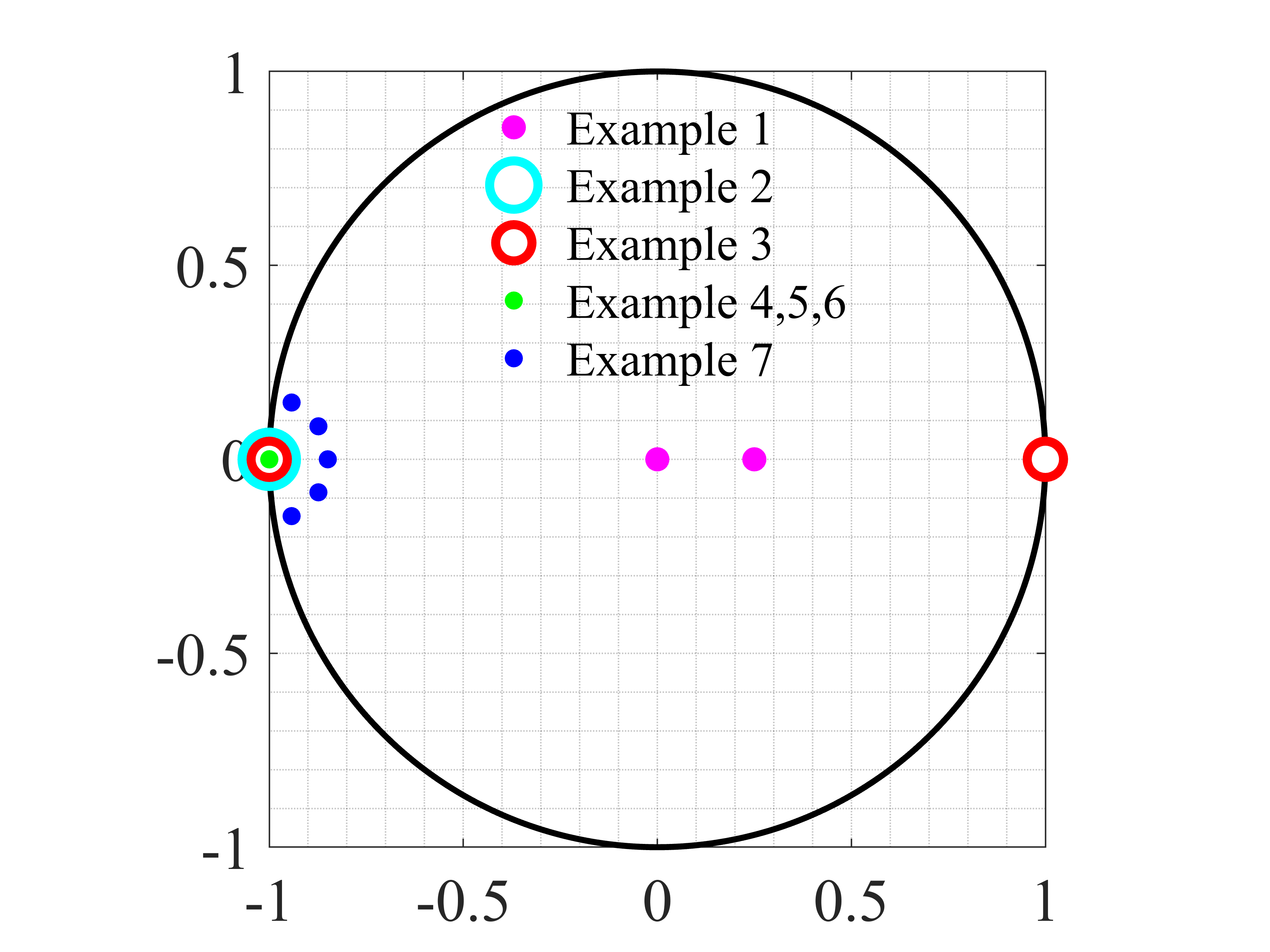}
  \caption{The eigenvalues of the matrix pencil for DARE for examples 1--7.}
  \label{fig:are1}
\end{figure}

As the results for Maple 17 in Table \ref{table:maple} and six different
Matlab versions in Table \ref{table:matlab} show, it is often impossible
to solve scalar and matrix Riccati equations using the built-in solvers
{\tt IDARE} or {\tt DARE} when the eigenvalues of the matrix pencil are
on the unit circle. The error using Maple is a a minimum of one order
larger than that of Matlab. Thus, the {\tt DARE} solver in Maple is not
recommended.

\begin{table}[cols=5,pos=ht]
  \caption{The errors $\epsilon_X$, $\epsilon_H$ and $\epsilon_P$ for
    examples 1--7, using the DARE or IDARE routine in several different
    Matlab
    versions. \\
    Note: The meaning of {\tt Report} has changed between versions. In
    older versions, {\tt Report=-1} meant ``the associated symplectic
    pencil has eigenvalues on or very near the unit circle''. In newer
    versions, {\tt Report=3} means ``The symplectic spectrum has
    eigenvalues on the unit circle''. These error messages are
    essentially the same.}\label{table:matlab}
\begin{tabular}{ccccc}
  \hline
  \hline
Matlab version & Examples & 	$\epsilon_X$ & $\epsilon_H$ & $\epsilon_P$ \\
  \hline
  \hline
  \multirow{10}{*}{2012a (DARE)}
               & 1 & 1.78e-15 & 5.55e-17 & 5.90e-16 \\
  \cline{2-5}
               & \multirow{4}{*}{2} & \multicolumn{3}{c}{$p_0 = 2$, $p_1 = 1$} \\
               & & \multicolumn{3}{c}{No solutions, Report = -1} \\
  \cline{3-5}
               & & \multicolumn{3}{c}{$p_0 = 1$, $p_1 = 1/2$} \\
               & & 4.95e-09 & 3.50e-09 & 0 \\
  \cline{2-5}
               & 3 & 1.49e-06 & 8.31e-07 & 5.68e-14 \\
               & 4 & \multicolumn{3}{c}{No solutions, Report = -1} \\
               & 5 & \multicolumn{3}{c}{No solutions, Report = -1}\\
               & 6 & \multicolumn{3}{c}{No solutions, Report = -1} \\
               & 7 & 0.00105 & 0.0137 & 1.09e-10 \\
  \hline
  \hline
  \multirow{10}{*}{R2018a (DARE)}
               & 1 & 7.11e-15 & 4.16e-17 & 1.78e-15 \\
  \cline{2-5}
               & \multirow{4}{*}{2} & \multicolumn{3}{c}{$p_0 = 2$, $p_1 = 1$} \\
               & & \multicolumn{3}{c}{No solutions, Report = -1} \\
  \cline{3-5}
               & & \multicolumn{3}{c}{$p_0 = 1$, $p_1 = 1/2$} \\
               & & 7.00e-09 & 4.95e-09 & 5.55e-17 \\
  \cline{2-5}
               & 3 & 1.56e-06 & 7.98e-07 & 9.95e-14 \\
               & 4 & 5.45e-05 & 5.45e-05 & 8.88e-16 \\
               & 5 & 3.58e-05 & 6.16e-05 & 6.66e-14 \\
               & 6 & 4.32e-05 & 6.12e-05 & 4.44e-16 \\
               & 7 & 0.00108 & 0.01418 & 1.978e-11 \\
  \hline
  \hline
%  \multirow{10}{*}{\mbox{R2020a \\ R2020b \\ R2021a}}
  \multirow{10}{*}{\mbox{\begin{tabular}{l} R2020a (IDARE) \\R2020b (IDARE) \\ R2021a (IDARE) \end{tabular}}}
               & 1 & 3.55e-15 & 1.11e-16 & 3.55e-15 \\
  \cline{2-5}
               & \multirow{4}{*}{2} & \multicolumn{3}{c}{$p_0 = 2$, $p_1 = 1$} \\
               & & \multicolumn{3}{c}{No solutions, Report = -1} \\
  \cline{3-5}
               & & \multicolumn{3}{c}{$p_0 = 1$, $p_1 = 1/2$} \\
               & & 5.55e-09 & 3.92e-09 & 0 \\
  \cline{2-5}
               & 3 & \multicolumn{3}{c}{No solutions, Report = 3} \\
               & 4 & \multicolumn{3}{c}{No solutions, Report = 3} \\
               & 5 & 4.45e-05 & 7.64e-05 & 2.291e-13 \\
               & 6 & 4.33e-05 & 6.12e-05 & 8.882e-16 \\
               & 7 & 0.00100 & 0.0131 & 1.462e-10 \\
  \hline
  \hline
%  \multirow{10}{*}{R2021a \\ Update 4} 
  \multirow{10}{*}{\mbox{\begin{tabular}{l} R2021a \\ Update 4 (IDARE) \end{tabular}}}
               & 1 & 3.55e-15 & 8.88e-16 & 3.55e-15 \\
  \cline{2-5}
               & \multirow{4}{*}{2} & \multicolumn{3}{c}{$p_0 = 2$, $p_1 = 1$} \\
               & & \multicolumn{3}{c}{No solutions, Report = 3} \\
  \cline{3-5}
               & & \multicolumn{3}{c}{$p_0 = 1$, $p_1 = 1/2$} \\
                & & \multicolumn{3}{c}{No solutions, Report = 3} \\
  \cline{2-5}
               & 3 & 8.54e-07 & 6.44e-07 & 1.56e-13 \\
               & 4 & \multicolumn{3}{c}{No solutions, Report = 3} \\
               & 5 & \multicolumn{3}{c}{No solutions, Report = 3}\\
               & 6 & 4.73e-05 & 6.69e-05 & 8.88e-16 \\
               & 7 & 0.00099 & 0.0128 & 2.97e-10 \\
  \hline
  \hline
  \multicolumn{5}{c}{\begin{minipage}{3in} 
  \end{minipage}}
\end{tabular}
\end{table}

Full precision is obtained only for the nonsingular case in example 1,
where the determinant of the second order polynomial matrix has zeros at
1/2 and 2. The precision of the spectral factors obtained for Example 2,
which has the double zeros on the unit circle, for all Matlab versions
is $\epsilon_H \approx 10^{-9}$. For the singular matrix examples, the
accurcay drops to $10^{-6}$ or $10^{-5}$, or even $10^{-3}$ in the
highly singular example 7. This leads to unacceptably low precision for
scaling and wavelet functions.
\end{reorganized}

\section{Conclusion}\label{sec:conclusion}

This paper provides a broad perspective of improvements on Bauer’s
method for matrix spectral factorization, especially in the singular
case. We introduced two fast numerical algorithms, and a symbolic
algorithm which produces exact factors for some low-order polynomial
matrices.

The fixed point iteration algorithm is equivalent to Bauer's method but
avoids forming large matrices.  The Newton's method algorithm with its
quadratic convergence is significantly faster.  The symbolic algorithm
allows exact solutions for the spectral density of a second order
stochastic process and low-order polynomial matrices in MIMO control
theory applications.

By applying matrix spectral factorization to the well–known CL product
filter we derive a new orthogonal multiscaling function with exact
coefficients, which leads to a new multiwavelet filter bank which
possesses better Sobolev smoothness and coding gain than the original CL
multiwavelet filter bank.  Based on its advantages, in a forthcoming
work the obtained multiwavelets will be considered in different
applications of numerical mathematics, and in the context of signal or
image processing.

As the results for different six Matlab versions in Table
\ref{table:matlab} show, it is often impossible to solve scalar and
matrix Riccati equations using the built-in solvers IDARE (DARE) when
the eigenvalues of the matrix pencil are on the unit circle.

\section{Acknowledgements}
The research of one of the authors was supported by the National Science
Foundation and the National Geospatial-Intelligence Agency under awards
\#DMS-1830254 and CCF-1750920.

The authors would like to thank the editor, two anonymous reviewers and
prof. Chun-Hua Guo, University of Regina, Canada, for their careful
reading, their constructive and valuable comments and very helpful
suggestions which helped us to improve the manuscript.
 
\appendix
\section*{Appendix A}
We have moved some technical proofs to this appendix, in order not to
interrupt the flow of arguments in the main paper.

\subsection*{A.1. The Block Cholesky Method}

It is well-known that matrix multiplication $A \cdot B = C$ can be
performed element by element, or equivalently block by block, if the
matrices have been subdivided into blocks of compatible size.

Consider the matrices $\calP$, $\calH$ from eqs.~\eqref{eq:P},
\eqref{eq:H}, with $\calP = \calH \cdot \calH^*$. We divide both
matrices into $m \times m$ blocks. For $m=3$, for example, we take the
element-by-element form
{\small
\begin{equation*}
  \begin{bmatrix}
    & & \ddots \\
    P_{-3} & P_{-2} & P_{-1} & P_0 & P_1 & P_2 & P_3 \\
    &  P_{-3} & P_{-2} & P_{-1} & P_0 & P_1 & P_2 & P_3 \\
    & &  P_{-3} & P_{-2} & P_{-1} & P_0 & P_1 & P_2 & P_3 \\
    & & & & & & \ddots
  \end{bmatrix} = 
  \begin{bmatrix}
    & \ddots \\
    H_3 & H_2 & H_1 & H_0 \\
    & H_3 & H_2 & H_1 & H_0 \\
    & & H_3 & H_2 & H_1 & H_0 \\
    & & & & \ddots
  \end{bmatrix}
  \begin{bmatrix}
    & H_3^* \\
    \ddots & H_2^* & H_3^* \\
    & H_1^* & H_2^* & H_3^* \\
    & H_0^* & H_1^* & H_2^* \\
    & & H_0^* & H_1^* & \ddots \\
    & & & H_0^*
  \end{bmatrix}
\end{equation*}
}
and write it in terms of $3 \times 3$ blocks:
{\small
\begin{equation*}
  \left[
    \begin{array}{ccccccccc}
      & & & \ddots \\
      \hline
      \multicolumn{1}{|c}{P_{-3}} & P_{-2} & P_{-1} &
      \multicolumn{1}{|c}{P_0} & P_1 & P_2 & 
      \multicolumn{1}{|c}{P_3} & & \multicolumn{1}{c|}{\ } \\
      \multicolumn{1}{|c}{\ } &  P_{-3} & P_{-2} & 
      \multicolumn{1}{|c}{P_{-1}} & P_0 & P_1 & 
      \multicolumn{1}{|c}{P_2} & P_3 & \multicolumn{1}{c|}{\ } \\
      \multicolumn{1}{|c}{\ } & &  P_{-3} & 
      \multicolumn{1}{|c}{P_{-2}} & P_{-1} & P_0 & 
      \multicolumn{1}{|c}{P_1} & P_2 & \multicolumn{1}{c|}{P_3} \\
      \hline
      & & & & & \ddots
    \end{array} \right] = 
  \left[
    \begin{array}{cccccc}
    & \ddots \\
    \hline
    \multicolumn{1}{|c}{H_3} & H_2 & H_1 & 
    \multicolumn{1}{|c}{H_0} & &  \multicolumn{1}{c|}{\ } \\
    \multicolumn{1}{|c}{\ } & H_3 & H_2 & 
    \multicolumn{1}{|c}{H_1} & H_0 & \multicolumn{1}{c|}{\ } \\
    \multicolumn{1}{|c}{\ } & & H_3 &
    \multicolumn{1}{|c}{H_2} & H_1 & \multicolumn{1}{c|}{H_0} \\
    \hline
    & & & & \ddots
    \end{array} \right] 
 \left[
   \begin{array}{ccccc}
     \cline{2-4}
     & \multicolumn{1}{|c}{H_3^*} & & \multicolumn{1}{c|}{\ } \\
    \ddots & \multicolumn{1}{|c}{H_2^*} & H_3^* & \multicolumn{1}{c|}{\ } \\
    & \multicolumn{1}{|c}{H_1^*} & H_2^* & \multicolumn{1}{c|}{H_3^*} \\[3pt]
    \cline{2-4}
    & \multicolumn{1}{|c}{H_0^*} & H_1^* & \multicolumn{1}{c|}{H_2^*} \\
    & \multicolumn{1}{|c}{\ } & H_0^* & \multicolumn{1}{c|}{H_1^*} & \ddots \\
     & \multicolumn{1}{|c}{\ } & & \multicolumn{1}{c|}{H_0^*} \\
     \cline{2-4}
  \end{array} \right].
\end{equation*}
}
This is equivalent to the block case with $m=1$
\begin{equation*}
  \hat\calP = \hat\calH \, \hat\calH^*
\end{equation*}
with
\begin{equation}\label{eq:blockPH}
  \hat\calP = 
  \begin{bmatrix}
    \ddots \\
    & \hat P_{-1} & \hat P_0 & \hat P_1 \\
    & & \hat P_{-1} & \hat P_0 & \hat P_1 \\
    & & & \ddots 
  \end{bmatrix}, \qquad
  \hat\calH = 
  \begin{bmatrix}
    \ddots \\
    & \hat H_1 & \hat H_0 \\
    & & \hat H_1 & \hat H_0 \\
    & & & \ddots 
  \end{bmatrix}, 
\end{equation}
where
\begin{equation*}
  \hat P_0 =
  \begin{bmatrix}
    P_0 & P_1 & P_2 \\[2pt]
    P_{-1} & P_0 & P_1 & \\[2pt]
    P_{-2} & P_{-1} & P_0 
  \end{bmatrix}, \quad
  \hat P_1 =
  \begin{bmatrix}
    P_3 & 0 & 0 \\
    P_2 & P_3 & 0 \\
    P_1 & P_2 & P_3
  \end{bmatrix}, \quad
  \hat H_0 =
  \begin{bmatrix}
    H_0 & 0  & 0  \\
    H_1 & H_0 & 0 & \\
    H_2 & H_1 & H_0 
  \end{bmatrix}, \quad
  \hat H_1 =
  \begin{bmatrix}
    H_3 & H_2 & H_1 \\
    0 & H_3 & H_2 \\
    0 & 0 & H_3
  \end{bmatrix}
\end{equation*}
and $\hat P_{-1} = \hat P_1^*$.

Next, we observe that the Cholesky decomposition in block form is
equivalent to the standard form.

As a simple example, we consider $2 \times 2$ case. Assume the matrix
\begin{equation}
  A =
  \begin{bmatrix}
    a_{11} & a_{12} \\
    a_{21} & a_{22}
  \end{bmatrix}
\end{equation}
is positive definite, so that $a_{11}$, $a_{22}$ are positive, and
$a_{12} = \overline{a_{21}}$.
In the Cholesky method, we set up
\begin{equation*}
  L =
  \begin{bmatrix}
    l_{11} & 0 \\
    l_{21} & l_{22}
  \end{bmatrix}
\end{equation*}
and demand
\begin{equation*}
  A = L  L^* =
  \begin{bmatrix}
    |l_{11}|^2 & l_{11} \overline{l_{21}} \\
     l_{21} \overline{l_{11}} &  |l_{11}|^2 + |l_{22}|^2
  \end{bmatrix}.
\end{equation*}
The entries in $L$ can then be computed one by one:
\begin{align*}
  l_{11} &= \sqrt{a_{11}} \\
  l_{21} &= a_{21}/ \overline{l_{11}} \\
  l_{22} &= \sqrt{ a_{22} -  |l_{21}|^2}
\end{align*}
Likewise, for a positive definite $2 \times 2$ block matrix
\begin{equation}
  A =
  \begin{bmatrix}
    A_{11} & A_{12} \\
    A_{21} & A_{22}
  \end{bmatrix} = L \, L^*, \qquad
    L =
  \begin{bmatrix}
    L_{11} & 0 \\
    L_{21} & L_{22}
  \end{bmatrix}
\end{equation}
we get
\begin{align*}
  L_{11} &= \mbox{Cholesky factor of $A_{11}$} \\
  L_{21} &= A_{21} L_{11}^{-*} \\
  L_{22} &= \mbox{Cholesky factor of $A_{22} - L_{21} L_{21}^*$}
\end{align*}
Cholesky factorization is unique as long as the factors are chosen to
have positive numbers on the diagonal. Since the element-by-element and
block versions of the algorithm both provide a Cholesky factorization
with positive diagonal, the results are identical.

This argument can be continued to larger matrices.

Let us now apply this to the factorization of
$\calP_n = \calH_n \calH_n^*$ (see eq.~\eqref{eq:Pn}). Beginning again
with the simplest case $n=2$
\begin{equation*}
  \calP_2 =
  \begin{bmatrix}
    P_0    & P_1 \\
    P_{-1} & P_0
  \end{bmatrix}, \qquad
  \calH_2 =
  \begin{bmatrix}
    H_0^{(0)}  & 0 \\[3pt]
    H_1^{(1)} & H_0^{(1)} 
  \end{bmatrix}.
\end{equation*}
We find
\begin{align*}
  H_0^{(0)} &= \mbox{Cholesky factor of $P_0$} \\
  H_1^{(1)} &= P_1^*  \left[ H_0^{(0)} \right]^{-*} \\
  H_0^{(1)} &= \mbox{Cholesky factor of $P_0 - H_1^{(1)} \left[ H_1^{(1)} \right]^{-*}$} \\
              &  =\mbox{Cholesky factor of $P_0 - P_1^* \left[ H_0^{(0) }\right]^{-*} \left[ H_0^{(0)} \right]^{-1} P_1$}.
\end{align*}
If we define
\begin{equation*}
  X^{(0)} = H_0^{(0)} \left[ H_0^{(0)} \right]^*, \qquad X^{(1)} = H_0^{(1)} \left[ H_0^{(1)} \right]^*,
\end{equation*}
this leads to
\begin{equation*}
  X^{(1)} = P_0 - P_1^* \left[ X^{(0)} \right]^{-1} P_1.
\end{equation*}
This process can be continued for larger matrices, to lead to
eq.~\eqref{eq:fpi}.

The fact that the last row in $\calH_n$ in Bauer's method converges to
the true solution immediately implies that the last $m$ rows
simultaneously converge, and the $H$ coefficients therefore
converge to the block form given in~ \eqref{eq:blockPH}.

\subsection*{A.2. Errors and Residuals}

We show in this subsection that in our particular case, it is enough to
consider the error in $H_0$ and the residual for $P_0$.

We assume that $H_0$, $H_1$ are the coefficients of the true factor
$H(z)$, and $H_0^{(n)}$, $H_1^{(n)}$ are the approximate factors at step
$n$:
\begin{gather*}
  H_0^{(n)} = H_0 + \Delta H_0^{(n)}, \\
  H_1^{(n)} = H_1 + \Delta H_1^{(n)}.
\end{gather*}
We define the error in $H$ as
\begin{equation*}
  \epsilon_H^{(n)} = \| \Delta H_0^{(n)} \| = \| H_0 - H_0^{(n)} \|
\end{equation*}
for some suitable matrix norm. It suffices to consider the error in
$H_0$, for the following reason.

Assuming that $A$ is an invertible matrix and $\Delta A$ is small
enough, we find
\begin{equation}\label{eq:ADAinv}
  \begin{split}
    (A + \Delta A)^{-1} &= \left[ A (I + A^{-1} \Delta A) \right]^{-1} = (I + A^{-1} \Delta A)^{-1} A^{-1} \\
      &\approx (I - A^{-1} \Delta A) \, A^{-1} = A^{-1} - A^{-1} \, \Delta A \, A^{-1}.
   \end{split}
\end{equation}
This gives us
\begin{equation*}
  [H_0^{(n)}]^{-*} = (H_0 + \Delta H_0^{(n)})^{-*} \approx H_0^{-*} -
  H_0^{-*} (\Delta H_0^{(n)})^* H_0^{-*},
\end{equation*}
so by \eqref{eq:findH1}
\begin{align*}
  H_1^{(n)} &= P_1^* [H_0^{(n)}]^{-*} \approx P_1^* \left[ H_0^{-*} -
    H_0^{-*} (\Delta H_0^{(n)})^* H_0^{-*} \right] = H_1 -  P_1^* H_0^{-*} (\Delta H_0^{(n)})^* H_0^{-*}, \\
  \Delta H_1^{(n)} &\approx -  P_1^* H_0^{-*} (\Delta H_0^{(n)})^* H_0^{-*}, \\
  \| \Delta H_1^{(n)}\| &\le  \|P_1\| \, \| H_0^{-1}\|^2 \, \| \Delta
  H_0^{(n)} \| = c \, \| \Delta H_0^{(n)} \|
\end{align*}
for a fixed constant $c$. Thus,
\begin{equation*}
  \| \Delta H_0^{(n)} \| + \| \Delta H_1^{(n)} \| \le (1+c) \| \Delta H_0^{(n)} \|.
\end{equation*}
Thus, the error in $H_0$ and the combined errors of $H_0$ and $H_1$
only differ by a constant factor.

We define likewise the residual
\begin{equation*}
  \epsilon_P^{(n)} = \| P_0 - P_0^{(n)} \|,
\end{equation*}
where $P_0^{(n)}$ is calculated from
\begin{equation*}
  P_0^{(n)} = H_0^{(n)} H_0^{(n)*} + H_1^{(n)} H_1^{(n)*}.
\end{equation*}
Since $H_1^{(n)}$ is computed from $P_1 = H_0^{(n)} [H_1^{(n)}]^{-*}$,
the residual of $P_1$ is zero.

It follows easily from \eqref{eq:taylor} that in the nonsingular case
$\epsilon_P^{(n)} = O(\epsilon_H^{(n)})$.

In the scalar singular case, we have
$\epsilon_P^{(n)} = O([\epsilon_H^{(n)}]^2)$. In our numerical
experiments, we observe that the same behavior still holds in the matrix
singular case.

\subsection*{A.3. Equivalent Concepts of Singularity}

In this appendix we show that the different concepts of singularity
listed in subsection \ref{subsec:singular} are identical.

We consider first the case of fixed point iteration, based on
\eqref{eq:g}.  Assume that the MSF problem is singular, that is,
$\left|P(\theta)\right| = 0$ for some $\theta \in \T$. Since
$P(\theta) = H(\theta) H(\theta)^*$, this implies that also
$\left|H(\theta)\right| = 0$. We find
\begin{equation*}
  \left| H(\theta) \right| = \left| H_0 + H_1 \thetabar \right|
  = \left| I + H_1 H_0^{-1} \thetabar \right| \, \left| H_0
  \right| = 0.
\end{equation*}
Since $H_0$ is nonsingular, the term $I + H_1 H_0^{-1} \thetabar$
must be singular. Thus, there exists a nonzero vector $\vecv$ which is
an eigenvector to eigenvalue 0 of $I + H_1 H_0^{-1} \thetabar$, or
equivalently an eigenvector to eigenvalue $(-\theta)$ of $H_1 H_0^{-1}$:
\begin{equation*}
  (I + H_1 H_0^{-1} \thetabar) \vecv = \veczero \quad
  \Leftrightarrow \quad H_1 H_0^{-1} \vecv = - \theta \vecv.
\end{equation*}

As in eq.~\eqref{eq:taylor2} we compute that
\begin{equation*}
  g(X + \Delta X) \approx g(X) + P_1^* X^{-1} (\Delta X) X^{-1} P_1.
\end{equation*}
After applying the $\vec$ operation and using \eqref{eq:vecABC}, we
get
\begin{equation*}
  \vecg(\vecx + \Delta \vecx) \approx \vecg(\vecx) + (P_1^T
  X^{-T})\otimes(P_1^* X^{-1}) \Delta\vecx.
\end{equation*}
Comparing this to the Taylor expansion
\begin{equation*}
  \vecg(\vecx + \Delta\vecx) \approx \vecg(\vecx) + D\,\vecg(\vecx) \Delta\vecx,
\end{equation*}
we see that
\begin{equation*}
  D\,\vecg(\vecx) = (X^{-T} P_1^T)\otimes(P_1^* X^{-1}).
\end{equation*}

At the solution $\overline X = H_0 H_0^*$ we have
\begin{equation*}
  P_1^* X^{-1} = \left(H_1 H_0^*\right)\,\left(H_0^{-*} H_0^{-1}\right) = H_1 H_0^{-1}.
\end{equation*}
Thus, the derivative matrix at the solution is
\begin{equation*}
  D\,\vecg(\vecxbar) = \overline{(H_1 H_0^{-1})} \otimes (H_1 H_0^{-1}),
\end{equation*}
which has eigenvalue $(-\thetabar)(-\theta) = 1$. The same
argument also works in reverse: if $D\,\vecg(\vecxbar)$ has an
eigenvalue of 1, then $P(\theta)$ is singular for some $\theta \in \T$.

This shows that the spectral radius $\rho(D\,\vecg(\vecxbar)) \ge 1$. To
complete the argument, we have to invoke the known convergence of the
algorithm from \cite{YK78}, which implies that the spectral radius
$\rho(D\,\vecg(\vecxbar)) \le 1$.

The statement for Newton's method is immediate: if
$\vecf(\vecx) = \vecx - \vecg(\vecx)$, then $D\,\vecg(\vecxbar)$ has an
eigenvalue 1 if and only if $D\,\vecf(\vecxbar)$ has an eigenvalue 0.

\bibliographystyle{cas-model2-names}
\bibliography{lo}
\end{document}